\documentclass[11pt]{article}
\usepackage{amsmath,amssymb,amsfonts}

\newtheorem{theorem}{Theorem}[section]
\newtheorem{proposition}[theorem]{Proposition}

\newtheorem{remark}[theorem]{Remark}
\newtheorem{corollary}[theorem]{Corollary}
\newtheorem{definition}[theorem]{Definition}
\newtheorem{example}[theorem]{Example}
\begin{document}

\title{Generalized geometric structures \\
on complex and symplectic manifolds}
\author{Marcos Salvai\thanks{%
Partially supported by Conicet, Foncyt, Secyt Univ.\thinspace Nac.\thinspace
C\'{o}rdoba}}
\date{}
\maketitle

\begin{abstract}
On a smooth manifold $M$, generalized complex (generalized paracomplex)
structures provide a notion of interpolation between complex (paracomplex)
and symplectic structures on $M$.

Given a complex manifold $\left( M,j\right) $, we define six families of
distinguished generalized complex or paracomplex structures on $M$. Each one
of them interpolates between two geometric structures on $M$ compatible with 
$j$, for instance, between totally real foliations and K\"{a}hler
structures, or between hypercomplex and $\mathbb{C}$-symplectic structures.
These structures on $M$ are sections of fiber bundles over $M$ with typical
fiber $G/H$ for some Lie groups $G$ and $H$. We determine $G$ and $H$ in
each case.

We proceed similarly for symplectic manifolds. We define six families of
generalized structures on $\left( M,\omega \right) $, each of them
interpolating between two structures compatible with $\omega $, for
instance, between a $\mathbb{C}$-symplectic and a para-K\"{a}hler structure
(aka bi-Lagrangian foliation).
\end{abstract}

\noindent 2000 MSC:\ 22F30, 22F50, 53B30, 53B35, 53C15, 53C56, 53D05

\smallskip

\noindent Keywords: generalized complex structure; interpolation; K\"{a}%
hler; hypercomplex; signature

\section{Introduction}

Generalized complex geometry arose from the work \cite{H} of Nigel Hitchin.
It has complex and symplectic geometry as its extremal special cases and
provides a notion of interpolation between them. It has greatly expanded
since its introduction only a decade ago and has far-reaching applications
in Mathematical Physics. We expect the development of similar ideas to be of
interest, now starting from a manifold which is already endowed with a
structure, and working out a notion of interpolation of supplementary
compatible geometric structures. Besides, we hope that, as it happens with
natural defined new structures, the search for nontrivial examples can
contribute, in same cases, to a better understanding of some manifolds, in
the same way, for instance, that generalized complex structures shed light
on the geometry of nil- and solvmanifolds \cite{Cav,Bar}.

Next we comment on the contents of the paper. In Section \ref{S2} we recall
the definitions and properties of generalized complex or paracomplex
structures. In Section \ref{S3} we have a manifold $M$ with a complex
structure $j$ and consider geometric structures on $M$ compatible with $j$,
which we call integrable $\left( \lambda ,0\right) $- or $\left( 0,\ell
\right) $-structures, with $\lambda ,\ell =\pm 1$; for instance, $\lambda
=-1 $ and $\ell =1$ give us hypercomplex and pseudo-K\"{a}hler structures,
respectively. The reason of this nomenclature is that it will allow us to
define families of generalized complex or paracomplex structures on $M$,
called integrable $\left( \lambda ,\ell \right) $-structures, which in a
certain sense, specified in Theorem \ref{THMcomplex}, interpolate between
integrable $\left( \lambda ,0\right) $- and $\left( 0,\ell \right) $%
-structures on $M$. In order to give strength to the notion of these
generalized structures on $M$, we prove that they are sections of fiber
bundles over $M$ with typical fiber $G/H$ for some Lie groups $G$ and $H$.
We determine $G$ and $H$ in each case. In Section \ref{S4} we proceed
similarly for a symplectic (instead of a complex) manifold.

\section{Generalized complex and paracomplex structures\label{S2}}

In this section, we recall from the seminal work \cite{gcg} the definitions
and basic facts on generalized complex structures, and on generalized
paracomplex structures (from \cite{Wade}).

Let $M$ be a smooth manifold (by smooth we mean of class $C^{\infty }$; all
the objects considered will belong to this class). The extended tangent
bundle is the vector bundle $\mathbb{T}M=TM\oplus TM^{\ast }$ over $M$. A
canonical split pseudo-Riemannian structure on $\mathbb{T}M$ is defined by 
\begin{equation*}
b\left( u+\sigma ,v+\tau \right) =\tau \left( u\right) +\sigma \left(
v\right) \text{,}
\end{equation*}%
for smooth sections $u+\sigma ,v+\tau $ of $\mathbb{T}M$. The \emph{Courant
bracket} of these sections \cite{Courant} is given by 
\begin{equation*}
\left[ u+\sigma ,v+\tau \right] =\left[ u,v\right] +\mathcal{L}_{u}\tau -%
\mathcal{L}_{v}\sigma -\tfrac{1}{2}d\left( \tau \left( u\right) -\sigma
\left( v\right) \right) \text{,}
\end{equation*}%
where $\mathcal{L}$ denotes the Lie derivative.

A \emph{paracomplex structure} $r$ on the smooth manifold $M$ is a smooth
tensor field of type $\left( 1,1\right) $ on $M$ satisfying $r^{2}=$ id such
that the eigendistributions of $r$ associated to the eigenvalues $1$ and $-1$
are integrable and have the same dimension \cite{Rocky}. Among all the
equivalent definitions of a \emph{complex structure} $j$ on $M$ we choose
the following:\ It is a smooth tensor field of type $\left( 1,1\right) $ on $%
M$ satisfying $j^{2}=-$ id such that the eigendistributions of $j$ in $%
TM\otimes \mathbb{C}$ associated to the eigenvalues $i$ and $-i$ are
involutive (for the $\mathbb{C}$-bilinear extension of the Lie bracket).

A real linear isomorphism $S$ with $S^{2}=\lambda $ id, $\lambda =\pm 1$, is
called \emph{split} if it has exactly two eigenspaces (of the
complexification of the vector space, if $\lambda =-1$) with the same
dimension; this is always the case if $\lambda =-1$.

For $\lambda =\pm 1$, let $S$ be a smooth section of End\thinspace $\left( 
\mathbb{T}M\right) $ satisfying 
\begin{equation}
S^{2}=\lambda \;\text{id, }S\text{ is split and skew-symmetric for }b
\label{condition1}
\end{equation}%
and such that the set of smooth sections of the $\pm \sqrt{\lambda }$%
-eigenspace of $S$ is closed under the Courant bracket (if $\lambda =-1$,
this means as usual closedness under the $\mathbb{C}$-linear extension of
the bracket to sections of the complexification of $\mathbb{T}M$). Then, for 
$\lambda =-1$ (respectively, $\lambda =1$), $S$ is called a \emph{%
generalized complex} (respectively, \emph{generalized paracomplex})
structure on $M$. Notice that in \cite{Wade} the latter is not required to
be split.

We also need the notion of $\left( +\right) $\emph{-generalized paracomplex}
structure $S$. It is the same as a generalized paracomplex structure, but
closedness under the Courant bracket is required only for sections of the
1-eingendistribution of $S$.

As far as we know, Izu Vaisman \cite{Izu}, was the first one to consider
generalized complex and paracomplex structures simultaneously in a
systematic way.

\section{Generalized geometric structures on complex manifolds\label{S3}}

\subsection{Geometric structures compatible with $j$}

Let $\left( M,j\right) $ be a complex manifold. We consider the following
well-known integrable geometric structures on $M$ compatible with $j$. The
reason of the names \emph{integrable }$\left( \lambda ,0\right) $\emph{- or }%
$\left( 0,\ell \right) $\emph{-structures} will become apparent in Theorem %
\ref{THMcomplex}.

\bigskip

\noindent \textbf{Integrable (}$1,0$\textbf{)-structure or complex product
structure on }$\left( M,j\right) $\textbf{.} It is given by a paracomplex
structure $r$ on $M$ with $rj=-jr$. Then $\left( M,j,r\right) $ is a complex
product manifold \cite{AdrianS}, also called \textbf{para-hypercomplex} \cite%
{Y, Davidov} or \textbf{neutral hypercomplex manifold} \cite{K, Dunaj}.

\medskip

\noindent \textbf{Integrable (}$-1,0$\textbf{)-structure or hypercomplex
structure on }$\left( M,j\right) $\textbf{.} It is given by a complex
structure $r$ on $M$ which is $j$-antilinear, that is, $rj=-jr$.

\medskip

\noindent \textbf{Integrable (}$0,1$\textbf{)-structure or pseudo-K\"{a}hler
structure on }$\left( M,j\right) $\textbf{.\ }It is given by a symplectic
form $\omega $ on $M$ for which $j$ is skew-symmetric. If $g$ denotes the
pseudo-Riemannian metric given by $g\left( u,v\right) =\omega \left(
ju,v\right) $, then $\left( M,g,j\right) $ is pseudo-K\"{a}hler with even
signature (since $j$ is an isometry for $g$).\textbf{\ }

\medskip

\noindent \textbf{Integrable (}$0,-1$\textbf{)-structure or }$\mathbb{C}$%
\textbf{-symplectic structure on }$\left( M,j\right) $\textbf{.} It is given
by a symplectic form $\omega $ on $M$ for which $j$ is symmetric. If $\theta 
$ denotes the two-form given by $\theta \left( u,v\right) =\omega \left(
ju,v\right) $, then $\Omega =\omega -i\theta $ is a $\mathbb{C}$-symplectic
structure on $M$.

\medskip

We also have

\medskip

\noindent \textbf{(}$\mathbf{+}$\textbf{)-integrable (}$1,0$\textbf{%
)-structure or totally real foliation of }$\left( M,j\right) $\textbf{.} It
is given by a tensor field $r$ of type $\left( 1,1\right) $ on $M$ with $%
r^{2}=\;$id and $rj=-jr$, such that the $1$-eigensection $\mathcal{D}$ of $r$
is an integrable distribution. Then $\mathcal{D}\oplus j\mathcal{D}=TM$
holds and the leaves of $\mathcal{D}$ are totally real submanifolds of $M$.

\bigskip

Recall that for a hypercomplex or a complex product structure $\left(
j,r\right) $, $jr$ turns out to be split and integrable (see \cite{AdrianS}%
). Also, if $j$ is a complex structure on $M$ and $\omega $ is a symplectic
form on $M$ for which $j$ is symmetric, then $\Omega $ (or equivalently $%
\theta $) as above is closed. Notice that hypercomplex and $\mathbb{C}$%
-symplectic manifolds have even complex dimension.

\subsection{Slash structures on $\left( M,j\right) $}

\begin{definition}
Let $\left( M,j\right) $ be a complex manifold. For $\ell =\pm 1$, let $%
J_{\ell }$ be the complex structure on the real vector bundle $\mathbb{T}M$
over $M$ given by 
\begin{equation*}
J_{\ell }=\left( 
\begin{array}{cc}
j & 0 \\ 
0 & \ell j^{*}%
\end{array}
\right) \text{.}
\end{equation*}
\end{definition}

Notice that $J_{-1}$ is a generalized complex structure on $M$, but $J_{1}$
is not, since\ it is not skew-symmetric for $b$. Indeed, for all sections $%
u+\sigma ,v+\tau $ of $\mathbb{T}M$, 
\begin{eqnarray}
b\left( J_{\ell }\left( u+\sigma \right) ,v+\tau \right) &=&b\left( ju+\ell
j^{\ast }\sigma ,v+\tau \right) =\tau \left( ju\right) +\ell \sigma \left(
jv\right)  \label{Jele} \\
&=&\ell b\left( u+\sigma ,J_{\ell }\left( v+\tau \right) \right) \text{.} 
\notag
\end{eqnarray}

\bigskip

Now we introduce six families of generalized geometric structures on $\left(
M,j\right) $ interpolating between some of the structures listed in the
previous subsection.

\begin{definition}
\label{Def}Let $\left( M,j\right) $ be a complex manifold. Given $\lambda
=\pm 1$ and $\ell =\pm 1$, a generalized complex structure $S$ \emph{(}for $%
\lambda =-1$\emph{)} or a generalized paracomplex structure $S$ \emph{(}for $%
\lambda =1$\emph{)} on $M$ is said to be an \textbf{integrable }$\left(
\lambda ,\ell \right) $\textbf{-structure} on $\left( M,j\right) $ if 
\begin{equation}
SJ_{\ell }=-J_{\ell }S\text{.}  \label{condition2}
\end{equation}%
Analogously, given $\ell =\pm 1$, a $\left( +\right) $-generalized
paracomplex structure $S$ on $M$ is said to be a $\left( +\right) $\textbf{%
-integrable }$\left( 1,\ell \right) $\textbf{-structure }on $\left(
M,j\right) $ if $SJ_{\ell }=-J_{\ell }S$.
\end{definition}

We call $\mathcal{S}_{j}\left( \lambda ,\ell \right) $ the set of all
integrable $\left( \lambda ,\ell \right) $-structures on $\left( M,j\right) $%
, and $\mathcal{S}_{j}^{+}\left( 1,\ell \right) $ the set of all $\left(
+\right) $-integrable $\left( 1,\ell \right) $-structures. An element of $%
\mathcal{S}_{j}\left( -1,-1\right) $ may be called, for instance, a
hypercomplex\thinspace /\thinspace $\mathbb{C}$-symplectic structure on $%
\left( M,j\right) $. That suggests the name \textbf{slash structures}\emph{\ 
}for these structures on $M$.

\medskip

Given a bilinear form $c$ on a real vector space $V$, let $c^{\flat }\in $
End\thinspace $\left( V,V^{\ast }\right) $ be defined by $c^{\flat }\left(
u\right) \left( v\right) =c\left( u,v\right) $. The form $c$ is symmetric
(respectively, skew-symmetric) if and only if $\left( c^{\flat }\right)
^{\ast }=c^{\flat }$ (respectively, $\left( c^{\flat }\right) ^{\ast
}=-c^{\flat }$).

\begin{example}
If $r$ and $\omega $ are integrable $\left( \lambda ,0\right) $- and $\left(
0,\ell \right) $-structures on $\left( M,j\right) $, respectively, then 
\begin{equation}
R=\left( 
\begin{array}{cc}
r & 0 \\ 
0 & -r^{\ast }%
\end{array}%
\right) \text{ \ \ \ and\ \ \ \ \ }Q=\left( 
\begin{array}{cc}
0 & \lambda \left( \omega ^{\flat }\right) ^{-1} \\ 
\omega ^{\flat } & 0%
\end{array}%
\right)  \label{trivialC}
\end{equation}

belong to $\mathcal{S}_{j}\left( \lambda ,\ell \right) $.
\end{example}

\medskip

The following simple theorem justifies the terminology introduced in the
section and includes the notion of interpolation. See comments on this
concept in subsection \ref{weak}.

\begin{theorem}
\label{THMcomplex}Let $\left( M,j\right) $ be a complex manifold. For $%
\lambda =\pm 1,\ell =\pm 1$, integrable $\left( \lambda ,\ell \right) $%
-structures on $\left( M,j\right) $ interpolate between integrable $\left(
\lambda ,0\right) $- and $\left( 0,\ell \right) $-structures on $\left(
M,j\right) $, that is, if 
\begin{equation*}
R=\left( 
\begin{array}{ll}
r & 0 \\ 
0 & t%
\end{array}%
\right) \text{ \ \ \ and\ \ \ \ }Q=\left( 
\begin{array}{cc}
0 & p \\ 
\omega ^{\flat } & 0%
\end{array}%
\right) 
\end{equation*}%
belong to $\mathcal{S}_{j}\left( \lambda ,\ell \right) $, then $r$ and $%
\omega $ are integrable $\left( \lambda ,0\right) $- and $\left( 0,\ell
\right) $-structures on $\left( M,j\right) $, respectively.

Also, for $\ell =\pm 1$, $\left( +\right) $-integrable $\left( 1,\ell
\right) $-structures interpolate between $\left( +\right) $-integrable $%
\left( 1,0\right) $- and integrable $\left( 0,\ell \right) $-structures on $%
\left( M,j\right) $.
\end{theorem}

\noindent \textbf{Proof. }We call $q=\omega ^{\flat }$. We know from (1.14),
(1.15) and Theorem 1.1 in \cite{Izu} (based on \cite{Crainic}) that if $R$ and $%
Q$ as above are both generalized complex (respectively, paracomplex)\
structures, then $\omega $ is a closed $2$-form and $r$ is a complex
structure (respectively, a tensor field of type $\left( 1,1\right) $ with $%
r^{2}=$ id and involutive eigendistributions) on $M$. Also, that $t=-r^{\ast
}$ and $p=-q^{-1}$ (respectively, $p=q^{-1}$).

Now, since $R$ and $Q$ anti-commute with $J_{\ell }$, one has that $rj=-jr$
and $qj=-\ell j^{\ast }q$. The first fact implies that $j$ interchanges the
eigendistributions of $r$ and so $r$ is a paracomplex structure ($r$ is
split). The second fact yields $q\left( ju\right) \left( v\right) =-\ell
q\left( u\right) \left( jv\right) $ for all vector fields $u,v$ on $M$.
Hence, $\omega \left( ju,v\right) =-\ell \omega \left( u,jv\right) $ for all 
$u,v$ and so $j$ is symmetric or skew-symmetric for $\omega $, depending on
whether $\ell =-1$ or $\ell =1$. Thus, $r$ and $\omega $ are integrable $%
\left( \lambda ,0\right) $- and $\left( 0,\ell \right) $-structures,
respectively.

Now assume that $\lambda =1$. Suppose that $R$ is a $\left( +\right) $%
-integrable $\left( 1,\ell \right) $-structure, that is, the $1$%
-eigensection 
\begin{equation*}
\mathcal{D}_{+}=\left\{ u+\sigma \mid ru=u,r^{\ast }\sigma =-\sigma \right\}
\end{equation*}
of $R$ is involutive with respect to the Courant bracket. In particular,
given vector fields $u$ and $v$ with $ru=u$ and $rv=v$, we have by
definition of the bracket that $r\left[ u,v\right] =\left[ u,v\right] $, and
hence the $1$-eigensection of $r$ is integrable (even if the $\left(
-1\right) $-eigensection of $R$ is not).

Finally, let $Q$ be as above, and for $\delta =\pm 1$ let 
\begin{equation*}
\mathcal{E}_{\delta }=\left\{ u+\delta \omega ^{\flat }u\mid u\in TM\right\}
\end{equation*}
be the $\delta $-eigensection of $Q$, and suppose that $\mathcal{E}_{+}$ is
Courant involutive. We see that then so is also $\mathcal{E}_{-}$. Indeed,
given two vector fields $u,v$ on $M$, by definition of the Courant bracket,
for $\delta =\pm 1$, there exists a $1$-form $\xi $ on $M$ such that 
\begin{equation*}
\left[ u+\delta \omega ^{\flat }u,v+\delta \omega ^{\flat }v\right] =\left[
u,v\right] +\delta \xi \text{.}
\end{equation*}%
The assertion follows, since $\xi =\omega ^{\flat }\left( \left[ u,v\right]
\right) $ if $\mathcal{E}_{+}$ is involutive. Consequently, $Q$ is a
generalized paracomplex structure and thus $\omega $ is an integrable $%
\left( 0,\ell \right) $-structure on $M$. \hfill\ $\square $

\begin{remark}
The choice of five compatible geometric structures on $\left( M,j\right) $
was strongly conditioned by Courant involutivity. For instance, we have not
considered anti-K\"{a}hler structures $g$ on $\left( M,j\right) $, i.e.
pseudo-Riemannian metrics $g$ for which $j$ is symmetric and parallel \emph{%
\cite{AK}}, since we have not been able to relate the integrability
condition \emph{(}that $j$ be parallel with respect to the Levi-Civita
connection of $g$\emph{)} to the Courant bracket.
\end{remark}

\subsection{A signature associated to integrable $\left( 1,1\right) $%
-structures on $\left( M,j\right) $}

\begin{proposition}
Let $S$ be an integrable $\left( 1,1\right) $-structure on a complex
manifold $\left( M,j\right) $ of complex dimension $m$. Then the form $\beta
_{S}$ on $\mathbb{T}M$ defined by $\beta _{S}\left( x,y\right) =b\left(
SJ_{+}x,y\right) $ is symmetric and has signature $\left( 2n,4m-2n\right) $
for some integer $n$ with $0\leq n\leq 2m$.
\end{proposition}

\noindent \textbf{Proof. }The form $\beta _{S}$ is symmetric since $S$ and $%
J_{+}$ anti-commute and are skew-symmetric and symmetric for $b$ (see (\ref%
{Jele})), respectively.

One has that $\left( SJ_{+}\right) ^{2}=$ id. For $\delta =\pm 1$, let $%
D_{\delta }$ be the $\delta $-eigensection of $SJ_{+}$. One verifies that $%
J_{+}\left( D_{+}\right) =D_{-}$, so $D_{+}$ and $D_{-}$ have both dimension 
$2m$.

For $\delta =\pm 1$ let $b^{\delta }:=\left. b\right| _{D_{\delta }\times
D_{\delta }}$ and $\beta ^{\delta }:=\left. \beta _{S}\right| _{D_{\delta
}\times D_{\delta }}$. One computes $b\left( D_{+},D_{-}\right) =0$; in
particular, by the orthogonality lemma (2.30 in \cite{Harvey}), $b^{\delta }$
is nondegenerate. Suppose that $b^{+}$ has signature $\left( n,2m-n\right) $%
. Hence $b^{-}$ has signature $\left( 2m-n,n\right) $ ($b$ is split). On the
other hand, one computes also that $b^{\delta }=\delta \beta ^{\delta }$.
Therefore the signature of $\beta _{S}$ is $\left( 2n,4m-2n\right) $, as
desired. \hfill \ $\square $

\begin{definition}
An integrable $\left( 1,1\right) $-structure $S$ on $\left( M,j\right) $ as
above is called an \textbf{integrable} $\left( 1,1;n\right) $\textbf{%
-structure}, and we write \emph{sig\thinspace }$\left( S\right) =n$. If $%
\beta _{S}$ is split \emph{(}or equivalently, $n=m$\emph{)}, by the next
proposition, the $\left( 1,1;n\right) $-structure is called a \emph{(}%
complex product\emph{)}/\emph{(}split K\"{a}hler\emph{)} structure on $%
\left( M,j\right) $.
\end{definition}

\begin{proposition}
\label{signC}\emph{a) }Let $r$ be an integrable $\left( 1,0\right) $%
-structure on $\left( M,j\right) $, that is, a complex product structure on $%
M$ compatible with $j$. Then 
\begin{equation*}
R=\left( 
\begin{array}{cc}
r & 0 \\ 
0 & -r^{*}%
\end{array}
\right)
\end{equation*}
is a $\left( 1,1;n\right) $-structure on $\left( M,j\right) $ if and only if 
$n=m$.

\smallskip

\emph{b) }Let $\omega $ be an integrable $\left( 0,1\right) $-structure on $%
\left( M,j\right) $. Then 
\begin{equation*}
Q=\left( 
\begin{array}{cc}
0 & \left( \omega ^{\flat }\right) ^{-1} \\ 
\omega ^{\flat } & 0%
\end{array}%
\right)
\end{equation*}%
is a $\left( 1,1;n\right) $-structure on $\left( M,j\right) $ if and only if
the pseudo-K\"{a}hler metric $g\left( u,v\right) =\omega \left( ju,v\right) $
on $M$ has signature $\left( n,2m-n\right) $. In particular, $n$ is even.
\end{proposition}

\noindent \textbf{Proof. }a) Since $rj=-jr$, we compute 
\begin{equation*}
\beta _{R}\left( u+\sigma ,v+\tau \right) =\tau \left( rju\right) +\sigma
\left( rjv\right) \text{.}
\end{equation*}
Now, $rj$ squares to the identity and is split (its $\left( -1\right) $- and 
$1$-eigensections are interchanged by $j$). Then, locally, there exists a
basis $\left\{ u_{1},\dots u_{2m}\right\} $ of $TM$ such that $rj\left(
u_{i}\right) =u_{i}$ for $1\leq i\leq m$ and $rj\left( u_{i}\right) =-u_{i}$
for $m<i\leq 2m$. Let $\left\{ \alpha _{1},\dots ,\alpha _{2m}\right\} $ be
the dual basis. Analyzing the signs of $\beta _{R}\left( u_{i}+\alpha
_{i},u_{i}+\alpha _{i}\right) $ and $\beta _{R}\left( u_{i}-\alpha
_{i},u_{i}-\alpha _{i}\right) $, one concludes that $\beta _{R}$ is split,
and this yields (a).

b) One computes 
\begin{equation*}
\beta _{Q}\left( u+\sigma ,v+\tau \right) =\omega \left( ju,v\right) +\tau
\left( (\omega ^{\flat })^{-1}j^{\ast }\sigma \right) =g\left( u,v\right)
+h\left( \sigma ,\tau \right) \text{,}
\end{equation*}%
where the symmetric form $h$ on $T^{\ast }M$ is defined by the last
equality. Now, 
\begin{eqnarray*}
\left( (\omega ^{\flat })^{\ast }h\right) \left( z,w\right) &=&h\left(
\omega ^{\flat }z,\omega ^{\flat }w\right) =\omega ^{\flat }\left( w\right)
\left( (\omega ^{\flat })^{-1}j^{\ast }\omega ^{\flat }\left( z\right)
\right) = \\
&=&-\omega ^{\flat }\left( w\right) \left( (\omega ^{\flat })^{-1}\omega
^{\flat }\left( jz\right) \right) =-\omega ^{\flat }\left( w\right) \left(
jz\right) = \\
&=&\omega \left( jz,w\right) =g\left( z,w\right) \text{,}
\end{eqnarray*}%
since for an integrable $\left( 0,1\right) $-structure $\omega $ on $\left(
M,j\right) $, $j$ is skew-symmetric for $\omega $, that is, $j^{\ast }\omega
^{\flat }=-\omega ^{\flat }j$. Therefore, if $\phi :TM\oplus TM\rightarrow 
\mathbb{T}M$ is defined by $\phi \left( u,z\right) =\left( u,\omega ^{\flat
}z\right) $, then 
\begin{equation*}
\phi ^{\ast }\beta _{Q}\left( \left( u,z\right) ,\left( v,w\right) \right)
=g\left( u,v\right) +g\left( z,w\right) \text{. }
\end{equation*}%
This implies the assertion of (b), since $\phi ^{\ast }\beta _{Q}$ and $%
\beta _{Q}$ have the same signature. \hfill\ $\square $

\subsection{The associated bundles over $\left( M,j\right) \label{SubS}$}

Let $\mathbb{L}$ denote the Lorentz numbers $a+\varepsilon b,$ $\varepsilon
^{2}=1$. Let $V$ be a vector space over $\mathbb{F}=\mathbb{R}$, $\mathbb{C}$%
, $\mathbb{L}$ or $\mathbb{H}$, where $\mathbb{H}=\mathbb{C}+\mathbf{j}%
\mathbb{C}$ are the quaternions (we consider right vector spaces over $%
\mathbb{H}$). Recall from \cite{Harvey} that an $\mathbb{R}$-bilinear map $%
C:V\times V\rightarrow \mathbb{F}$ satisfying $C\left( x\lambda ,y\mu
\right) =\overline{\lambda }C\left( x,y\right) \mu $ for any $\lambda ,\mu
\in \mathbb{F}$ and $x,y\in V$ is called Hermitian (respectively,
anti-Hermitian) if $\overline{C\left( x,y\right) }=C\left( y,x\right) $
(respectively, $\overline{C\left( x,y\right) }=-C\left( y,x\right) $) for
all $x,y\in V$. Also, a Hermitian form on a vector space $V$ over $\mathbb{F}%
\ne \mathbb{L}$ is said to be split if it has $\mathbb{F}$-signature $\left(
n,n\right) $, where $2n=\dim _{\mathbb{F}}V$. The $\mathbb{L}$-signature
does not make sense, since $\overline{\varepsilon }\varepsilon =-1$.

\smallskip

Generalized complex structures on a $\left( 2n\right) $-dimensional manifold 
$N$ are sections of a bundle over $N$ with typical fiber $O\left(
2n,2n\right) /U\left( n,n\right) $ \cite{gcg}. In the same way, generalized
paracomplex structures on an $m$-dimensional manifold $N$ are sections of a
bundle over $N$ with typical fiber $O\left( m,m\right) /Gl\left( m,\mathbb{R}%
\right) $, since $Gl\left( m,\mathbb{R}\right) $ is the $\mathbb{L}$-unitary
group (Section 3 in \cite{HarveyL}). Theorem \ref{complexFibre} below
presents analogous statements for integrable $\left( \lambda ,\ell \right) $%
-structures on a complex manifold $\left( M,j\right) $.

\smallskip

Let $O\left( m,m\right) $ and $Sp\left( m,\mathbb{R}\right) $ be the groups
of automorphisms of a split symmetric and skew-symmetric form on $\mathbb{R}%
^{2m}$, respectively. Let $SO^{*}\left( 2m\right) $ and $Sp\left( n,n\right) 
$ ($m=2n)$ be the groups of automorphisms of an anti-Hermitian
(respectively, a split Hermitian) form on $\mathbb{H}^{m}$. In \cite{Harvey}
they are called $SK\left( m,\mathbb{H}\right) $ and $HU\left( n,n\right) $,
respectively.

\begin{theorem}
\label{complexFibre}Let $\left( M,j\right) $ be a complex manifold of
complex dimension $m$. Then, integrable $\left( \lambda ,\ell \right) $- or $%
\left( 1,1;n\right) $-structures on $\left( M,j\right) $ are smooth sections
of a fiber bundle over $M$ with typical fiber $G/H$, according to the
following table \emph{(}$m=2k$ in the case $\lambda =\ell =-1$\emph{)}. 
\begin{equation*}
\begin{tabular}{|c|c|c|c|c|}
\hline
$\lambda $ & $\ell $ & \emph{sig} & $G$ & $H$ \\ \hline
$1$ & $1$ & $n$ & $O\left( 2m,\mathbb{C}\right) $ & $O\left( n,2m-n\right) $
\\ \hline
$1$ & $-1$ & - & $U\left( m,m\right) $ & $Sp\left( m,\mathbb{R}\right) $ \\ 
\hline
$-1$ & $1$ & - & $O\left( 2m,\mathbb{C}\right) $ & $SO^{*}\left( 2m\right) $
\\ \hline
$-1$ & $-1$ & - & $U\left( 2k,2k\right) $ & $Sp\left( k,k\right) $ \\ \hline
\end{tabular}%
\end{equation*}
\end{theorem}

\begin{corollary}
A complex manifold admitting a hypercomplex~/$~\mathbb{C}$-symplectic
structure has even complex dimension.
\end{corollary}

Before proving the theorem we introduce some notation and present a
proposition. Now we work at the algebraic level. We fix $p\in M$ and call $%
\mathbb{E}=\mathbb{T}_{p}M$. By abuse of notation, in the rest of the
subsection we write $b$ and $J_{\ell }$ instead of $b_{p}$ and $\left(
J_{\ell }\right) _{p}$, omitting the subindex $p$. Also, we sometimes
identify $\left( 1,-1\right) =\left( +,-\right) $, etc.

Let $\sigma \left( \lambda ,\ell \right) $ denote the set of all $S\in $ End$%
\,_{\mathbb{R}}\left( \mathbb{E}\right) $ satisfying%
\begin{equation*}
S^{2}=\lambda \,\text{id, }S\text{ is split, skew-symmetric for }b\text{ and 
}SJ_{\ell }=-J_{\ell }S\text{.}
\end{equation*}

Note that $\left( \mathbb{E},J_{\ell }\right) $ is a vector space over $%
\mathbb{C}$ via $\left( a+ib\right) x=ax+J_{\ell }x$.

\begin{proposition}
\label{sesqui}For $\ell =\pm 1$, let $b_{\ell }:\mathbb{E}\times \mathbb{E}%
\rightarrow \mathbb{C}$ be defined by 
\begin{equation*}
b_{\ell }\left( x,y\right) =b\left( x,y\right) -ib\left( x,J_{\ell }y\right) 
\text{.}
\end{equation*}
Then $b_{-}$ is split $\mathbb{C}$-Hermitian and $b_{+}$ is $\mathbb{C}$%
-bilinear symmetric \emph{(}with respect to $J_{-},J_{+}$, respectively\emph{%
)}.

Also, if $S\in $ \emph{End}\thinspace $_{\mathbb{R}}\left( \mathbb{E}\right) 
$ satisfies $S^{2}=\lambda $ \emph{id}, then $S\in \sigma \left( \lambda
,\ell \right) $ if and only if 
\begin{equation}
b_{\ell }\left( Sx,Sy\right) =-\lambda \overline{b_{\ell }\left( x,y\right) }
\label{beleBarra}
\end{equation}
for any $x,y\in \mathbb{E}$.
\end{proposition}

\smallskip

\noindent \textbf{Proof. }First notice that $T\in $ End$\,_{\mathbb{R}%
}\left( \mathbb{E}\right) $ with $T^{2}=\mu ~$id is symmetric or
skew-symmetric for $b$ if and only if 
\begin{equation}
b\left( Tx,Ty\right) =\pm b\left( x,T^{2}y\right) =\pm \mu b\left( x,y\right)
\label{symm}
\end{equation}%
for all $x,y$. Using (\ref{Jele}) together with (\ref{symm}) with $T=J_{\ell
}$ and $\mu =\ell $, it is easy to check that 
\begin{equation*}
ib_{\ell }\left( x,y\right) =b_{\ell }\left( x,J_{\ell }y\right) =\ell
b_{\ell }\left( J_{\ell }x,y\right)
\end{equation*}%
for all $x,y$. Also, it follows immediately from the definitions that $%
b_{\ell }\left( x,y\right) =b_{\ell }\left( y,x\right) $ or $\overline{%
b_{\ell }\left( x,y\right) }=b_{\ell }\left( y,x\right) $ for all $x,y$,
depending on whether $\ell =1$ or $\ell =-1$, respectively. Besides, $b_{-}$
is split since $b=$ Re~$b_{-}$ is split. Thus, the first assertion is true.

Now we prove the second assertion. Suppose first that $S\in \sigma \left(
\lambda ,\ell \right) $. Since $S$ anti-commutes with $J_{\ell }$, we
compute (using (\ref{symm}) with $T=S$ and $\mu =\lambda $) 
\begin{eqnarray*}
b_{\ell }\left( Sx,Sy\right) &=&b\left( Sx,Sy\right) -ib\left( Sx,J_{\ell
}Sy\right) \\
&=&-\lambda b\left( x,y\right) +ib\left( Sx,SJ_{\ell }y\right) \\
&=&-\lambda b\left( x,y\right) -\lambda ib\left( x,J_{\ell }y\right) \\
&=&-\lambda \left( b\left( x,y\right) +ib\left( x,J_{\ell }y\right) \right)
\\
&=&-\lambda \overline{b_{\ell }\left( x,y\right) }\text{.}
\end{eqnarray*}%
Conversely, suppose that $S^{2}=\lambda $ id and (\ref{beleBarra}) holds. By
(\ref{symm}) with $T=S$ and $\mu =\lambda $, $S$ is skew-symmetric for $b=$
Re~$b_{\ell }$. Now we compute 
\begin{eqnarray*}
b_{\ell }\left( x,SJ_{\ell }y\right) &=&\lambda b_{\ell }\left(
S^{2}x,SJ_{\ell }y\right) =\lambda \left( -\lambda \right) \overline{b_{\ell
}\left( Sx,J_{\ell }y\right) }=-\overline{ib_{\ell }\left( Sx,y\right) }= \\
&=&-\left( -i\right) \lambda \overline{b_{\ell }\left( Sx,S^{2}y\right) }%
=i\lambda \left( -\lambda \right) b_{\ell }\left( x,Sy\right) =-b_{\ell
}\left( x,J_{\ell }Sy\right) \text{. }
\end{eqnarray*}%
Since $b_{\ell }$ is nondegenerate, $S$ anti-commutes with $J_{\ell }$. This
implies, in particular, that if $\lambda =1$, then $J_{\ell }\left(
D_{+}\right) =D_{-}$, where $D_{\pm }$ is the $\left( \pm 1\right) $%
-eigenspace of $S$. Hence, $S$ is split. Therefore, $S\in \sigma \left(
\lambda ,\ell \right) $. \hfill\ $\square $

\bigskip

The core of the arguments in the proofs of Theorems \ref{complexFibre} and %
\ref{symplFibre} is essentially from 1.6 in \cite{neretin}, except those
involving the Lorentz numbers. We put them in context and complete details
(write in coordinates, choose particular presentations, prove the
transitivity of the actions).

We use the notation and the standard forms of inner products of the book 
\cite{Harvey}. In particular, Hermitian and anti-Hermitian forms differ from
those in \cite{neretin} by conjugation. We resort repeatedly to the Basis
Theorem (\cite{Harvey}, 4.2). For inner products on $\mathbb{L}$-vector
spaces we refer to \cite{HarveyL} (where Lorentz numbers are called double
numbers and denoted by $\mathbb{D}$).

\bigskip

\noindent \textbf{Proof of Theorem \ref{complexFibre}. }For $\ell =\pm 1$,
by the first assertion in Proposition \ref{sesqui} and the Basis Theorem,
there exist complex linear coordinates $\phi _{\ell }^{-1}=\left( z,w\right)
:\left( \mathbb{E},J_{\ell }\right) \rightarrow \mathbb{C}^{2m}$ such that $%
B_{\ell }:=\phi _{\ell }^{*}b_{\ell }$ have the forms 
\begin{equation*}
B_{-}\left( \left( z,w\right) ,\left( z^{\prime },w^{\prime }\right) \right)
=\overline{z}^{t}z^{\prime }-\overline{w}^{t}w^{\prime }\ \ \ \ \ \text{and\
\ \ \ \ }B_{+}\left( Z,Z^{\prime }\right) =Z^{t}Z^{\prime }\text{,}
\end{equation*}
where $z,w,z^{\prime },w^{\prime }\in \mathbb{C}^{m},Z,Z^{\prime }\in 
\mathbb{C}^{2m}$ are column vectors and the superscript $t$ denotes
transpose.

Let $\Sigma \left( \lambda ,\ell \right) $ be the subset of End$\,_{\mathbb{R%
}}\left( \mathbb{C}^{2m}\right) $ corresponding to $\sigma \left( \lambda
,\ell \right) $ via the isomorphism $\phi _{\ell }$. By the second statement
of Proposition \ref{sesqui}, $U\left( m,m\right) $ and $O\left( 2m,\mathbb{C}%
\right) $ (the Lie groups preserving $B_{-}$ and $B_{+}$, respectively) act
by conjugation on $\Sigma \left( \lambda ,-\right) ,\Sigma \left( \lambda
,+\right) $, respectively.

In what follows, for each case $\left( \lambda ,\ell \right) \neq \left(
1,1\right) $ we present a particular real isomorphism $S$ of $\mathbb{C}%
^{2m} $ and show, using the second statement of Proposition \ref{sesqui},
that $S$ belongs to $\Sigma \left( \lambda ,\ell \right) $ (actually, we
write down the computation only for $\lambda =1=-\ell $, the other being
analogous). Then we check that the group $G$ associated to $\left( \lambda
,\ell \right) $ in the table acts transitively on $\Sigma \left( \lambda
,\ell \right) $, with isotropy subgroup the corresponding group $H$ in the
table. In this way, one concludes that $\Sigma \left( \lambda ,\ell \right) $
may be identified with $G/H$, as desired. The case $\left( 1,1;n\right) $ is
dealt with similarly.

\medskip

\textbf{Case (}$+,-$\textbf{):} Let $S\in $ End$\,_{\mathbb{R}}$\thinspace $%
\left( \mathbb{C}^{2m}\right) $ be defined by $S\left( z,w\right) =\left( 
\overline{w},\overline{z}\right) $.\ We use the second statement of
Proposition \ref{sesqui} to show that $S$ belongs to $\Sigma \left(
+,-\right) $. Clearly, $S^{2}=$ id and also 
\begin{equation*}
B_{-}\left( S\left( z,w\right) ,S\left( z^{\prime },w^{\prime }\right)
\right) =B_{-}\left( \left( \overline{w},\overline{z}\right) ,\left( 
\overline{w^{\prime }},\overline{z^{\prime }}\right) \right) =w^{t}\overline{%
w^{\prime }}-z^{t}\overline{z^{\prime }}=-\overline{B_{-}\left( \left(
z,w\right) ,\left( z^{\prime },w^{\prime }\right) \right) }\text{.}
\end{equation*}%
Now let $V$ be the $1$-eigenspace of $S$, that is, $V=\left\{ \left( z,%
\overline{z}\right) \mid z\in \mathbb{C}^{m}\right\} \cong \mathbb{R}^{2m}$.
One has $V\oplus iV=\mathbb{C}^{2m}$ and verifies that $\alpha :=-i\left.
B_{-}\right\vert _{V\times V}$ is a symplectic form on $V$. Indeed, one
computes $\alpha \left( \left( z,\overline{z}\right) ,\left( z^{\prime },%
\overline{z^{\prime }}\right) \right) =2\left( x^{t}y^{\prime
}-y^{t}x^{\prime }\right) $ if $z=x+iy$ and $z^{\prime }=x^{\prime
}+iy^{\prime }$.

Given $A\in Sp\left( V,\alpha \right) $, the map $\tilde{A}$ defined by $%
\tilde{A}\left( X+iY\right) =AX+iAY$, for $X,Y\in V$, is in $U\left(
m,m\right) $. This gives an inclusion of $Sp\left( m,\mathbb{R}\right) \cong
Sp\left( V,\alpha \right) $ into $U\left( m,m\right) $.

Now we check that the isotropy subgroup $H$ at $S$ of the action of $U\left(
m,m\right) $ on $\Sigma \left( +,-\right) $ is $Sp\left( V,\alpha \right) $.
Assume that $A\in Sp\left( V,\alpha \right) $. Clearly, $AS=SA$ ($\left.
S\right\vert _{V}=$ id$_{V}$). Hence, $\tilde{A}$ commutes with $S$, since $%
S $ is anti-linear. Then $\tilde{A}\in H$. Conversely, if $L\in U\left(
m,m\right) $ commutes with $S$, then $L$ preserves $V$ and so $L=\tilde{A}$
for some $A\in Sp\left( V,\alpha \right) $.

It remains to show that the action is transitive. Let $T\in \Sigma \left(
+,-\right) $ and let $W\ $be the $1$-eigenspace of $T$. One verifies, using (%
\ref{beleBarra}), that $\theta =-i\left. B_{-}\right| _{W\times W}$ is a
symplectic form on $W$. Let $X_{1},\dots ,X_{m},Y_{1},\dots ,Y_{m}$ be
vectors in $W$ such that $\theta \left( X_{s},Y_{t}\right) =2\delta _{st}$
and $\theta \left( X_{s},X_{t}\right) =\theta \left( Y_{s},Y_{t}\right) =0$
for all $s\leq t$. Let $F:V\rightarrow W$ be the linear transformation with $%
F\left( e_{s},e_{s}\right) =X_{s}$, $F\left( ie_{t},-ie_{t}\right) =Y_{t}$,
where $\left\{ e_{1},\dots ,e_{m}\right\} $ is the canonical basis of $%
\mathbb{R}^{m}$. Then $F$ extends $\mathbb{C}$-linearly to $\tilde{F}\in
U\left( m,m\right) $ such that $T=\tilde{F}S\tilde{F}^{-1}$. Therefore, $%
\Sigma \left( +,-\right) $ can be identified with $U\left( m,m\right)
/Sp\left( m,\mathbb{R}\right) $, as desired.

\medskip

\textbf{Case (}$-,-$\textbf{):} Any $S\in \Sigma \left( -,-\right) $ gives $%
\mathbb{C}^{2m}$ the structure of a right $\mathbb{H}$-vector space via $%
Z\left( u+\mathbf{j}v\right) =uZ+v\left( SZ\right) $ ($Z\in \mathbb{C}%
^{2m},u,v\in \mathbb{C}$). Given $S\in \Sigma \left( -,-\right) $, let 
\begin{equation*}
C\left( Z,Z^{\prime }\right) =B_{-}\left( Z,Z^{\prime }\right) -B_{-}\left(
Z,SZ^{\prime }\right) \mathbf{j}\text{.}
\end{equation*}%
By Lemma 2.72 in \cite{Harvey} (using (\ref{beleBarra}) and the fact that $u%
\mathbf{j}=\mathbf{j}\overline{u}$ for all $u\in \mathbb{C}$), $C$ is an $%
\mathbb{H}$-Hermitian form, which is split since $B_{-}$ is so. In
particular $m$ is even, say, $m=2k$. Now, $L\in U\left( m,m\right) $
commutes with $S$ if and only if $L$ is an isometry for $C$. Hence, the
isotropy subgroup at $S$ of the action of $U\left( m,m\right) $ is $Sp\left(
k,k\right) $. The action is transitive: If $T$ is another element of $\Sigma
\left( -,-\right) $, then one has another $\mathbb{H}$-structure on $\mathbb{%
E}$ and can define $C_{T}$ in the same way as $C$. By the Basis Theorem,
they are isometric. There exists an $\mathbb{H}$-linear isometry $F:\left( 
\mathbb{E},C\right) \rightarrow \left( \mathbb{E},C_{T}\right) $, which
satisfies $F\in U\left( m,m\right) $ and $T=FSF^{-1}$. Therefore, $\Sigma
\left( -,-\right) $ can be identified with $U\left( m,m\right) /Sp\left(
k,k\right) $, as desired.

We give an example of $S\in \Sigma \left( -,-\right) $:\ Write $z=\left(
z_{1},z_{2}\right) ,w=\left( w_{1},w_{2}\right) $, with $z_{s},w_{t}\in 
\mathbb{C}^{k}$ and define $S\in $ End$\,_{\mathbb{R}}\left( \mathbb{C}%
^{4k}\right) $ by $S\left( z_{1},z_{2},w_{1},w_{2}\right) =\left( -\overline{%
z_{2}},\overline{z_{1}},-\overline{w_{2}},\overline{w_{1}}\right) $.

\medskip

\textbf{Case (}$+,+,n$\textbf{): }Let $S\in $ End$\,_{\mathbb{R}}\left( 
\mathbb{C}^{2m}\right) $ be defined by $S\left( z,w\right) =\left( i%
\overline{z},-i\overline{w}\right) $, for $z\in \mathbb{C}^{n},w\in \mathbb{C%
}^{2m-n}$, which belongs to $\Sigma \left( +,+;n\right) $. In fact, one uses
the second statement of Proposition \ref{sesqui} to show that $S\in \Sigma
\left( +,+\right) $ and computes 
\begin{equation*}
\text{Re~}B_{+}\left( S\left( iz,iw\right) ,\left( z^{\prime },w^{\prime
}\right) \right) =\text{ Re~}\left( \overline{z}^{t}z^{\prime }-\overline{w}%
^{t}w^{\prime }\right) \text{,}
\end{equation*}%
which is a real symmetric form of signature $\left( 2n,4m\right) $. This
implies that $S\in \Sigma \left( +,+;n\right) $, since $b=$ Re\thinspace $%
b_{+}$. Let $V$ be the $1$-eigenspace of $S$, that is, 
\begin{equation*}
V=\left\{ \left( \left( 1+i\right) x,\left( 1-i\right) y\right) \mid x\in 
\mathbb{R}^{n},y\in \mathbb{R}^{2m-n}\right\} \cong \mathbb{R}^{2m}.
\end{equation*}%
Then $V\oplus iV=\mathbb{C}^{2m}$. One verifies that $g:=-i\left.
B_{+}\right\vert _{V\times V}$ is a real symmetric form on $V$ of signature $%
\left( n,2m-n\right) $. Indeed, one computes 
\begin{equation*}
g\left( \left( \left( 1+i\right) x,\left( 1-i\right) y\right) ,\left( \left(
1+i\right) x^{\prime },\left( 1-i\right) y^{\prime }\right) \right) =2\left(
x^{t}x^{\prime }-y^{t}y^{\prime }\right) \text{.}
\end{equation*}%
Given $A\in O\left( V,g\right) $, then $\tilde{A}\left( X+iY\right) =AX+iAY$
($X,Y\in V$) satisfies $\tilde{A}\in O\left( 2m,\mathbb{C}\right) $. This
gives an inclusion of $O\left( n,2m-n\right) $ in $O\left( 2m,\mathbb{C}%
\right) $.

We check that the isotropy subgroup at $S$ is $O\left( V,g\right) $: Since $%
S $ is anti-linear, $S$ commutes with $\tilde{A}$ for any $A\in O\left(
V,g\right) $. Besides, if $L\in O\left( 2m,\mathbb{C}\right) $ commutes with 
$S$, then $L$ preserves $V$ and so $L=\tilde{A}$ for some $A\in O\left(
V,g\right) $.

Now we see that the action is transitive. Let $T\in \Sigma \left(
+,+,n\right) $ and let $W$ be the $1$-eigenspace of $T$. Then $h:=-i\left.
B_{+}\right| _{W\times W}$ is a real symmetric form on $W$ of signature $%
\left( n,2m-n\right) $. In fact, if $Tu=u$ and $Tv=v$, 
\begin{equation*}
h\left( u,v\right) =-iB_{+}\left( u,v\right) =-iB_{+}\left( Tu,Tv\right) =i%
\overline{B_{+}\left( u,v\right) }=\overline{h\left( u,v\right) }
\end{equation*}
and also Re\thinspace $B_{+}\left( Tiu,v\right) =$ Re\thinspace $B_{+}\left(
-iTu,Tv\right) =h\left( u,v\right) $. Let $v_{1},\dots ,v_{2m}$ be a basis
of $W$ such that $h\left( v_{s},v_{s}\right) =2$ if $s\leq n$, $h\left(
v_{s},v_{s}\right) =-2$ if $s>n$ and $h\left( v_{s},v_{t}\right) =0$ for all 
$s\neq t$. Let $F:V\rightarrow W$ be the linear transformation with $F\left(
\left( 1+i\right) e_{s}\right) =v_{s}$ if $s\leq n$ and $F\left( \left(
1-i\right) e_{s}\right) =v_{s}$ if $s>n$. Then $F$ extends linearly to $%
\tilde{F}\in O\left( 2m,\mathbb{C}\right) $ such that $T=\tilde{F}S\tilde{F}%
^{-1}$.

\smallskip

\textbf{Case (}$-,+$\textbf{): }Let $S\in $ End$\,_{\mathbb{R}}\left( 
\mathbb{C}^{2m}\right) $ be defined by $S\left( z,w\right) =\left( -%
\overline{w},\overline{z}\right) $, which belongs to $\Sigma \left(
-,+\right) $. Notice that $\left( \mathbb{C}^{2m},S\right) $ is a right $%
\mathbb{H}$-vector space via $Z\left( z+\mathbf{j}w\right) =Zz+\left(
SZ\right) w$.

Let $C\left( Z,Z^{\prime }\right) =B_{+}\left( SZ,Z^{\prime }\right) -%
\mathbf{j}B_{+}\left( Z,Z^{\prime }\right) $. Then $C$ is skew $\mathbb{H}$%
-Hermitian. Now, $L\in O\left( 2m,\mathbb{C}\right) $ commutes with $S$ if
and only if $L$ is an isometry for $C$. Hence, the isotropy subgroup at $S$
of the action of $O\left( 2m,\mathbb{C}\right) $ is $SO^{*}\left( 2m\right) $%
. The action is transitive: If $T$ is another element of $\Sigma \left(
-,+\right) $, then one has another $\mathbb{H}$-structure on $\mathbb{E}$
compatible with $j$ and can define $C_{T}$ in the same way as $C$. By the
Basis Theorem, they are isometric. Then there exists an $\mathbb{H}$-linear
isometry $F:\left( \mathbb{E},C\right) \rightarrow \left( \mathbb{E}%
,C_{T}\right) $ satisfying $F\in O\left( 2m,\mathbb{C}\right) $ and $%
T=FSF^{-1}$. \hfill \ $\square $

\subsection{Slash structures and the notion of interpolation\label{weak}}

Generalized complex geometry on smooth manifolds generalizes complex and
symplectic structures and simultaneously provides a notion of interpolation
between them.

In our opinion, integrable $\left( \lambda ,\ell \right) $-structures on
complex manifolds are good generalizations of integrable $\left( \lambda
,0\right) $- or $\left( 0,\ell \right) $-structures, but for the sake of
simplicity, we have presented a rather weak definition of interpolation,
which in some cases is not what one would expect from that concept, but
(again in our view) in most cases is appropriate.

In the papers devoted to generalized complex structures the notion of
interpolation is not made explicit; there is no need of doing so, because
their existence on a smooth manifold $M$ implies the existence of almost
complex and almost symplectic structures on $M$. In contrast, on a complex
manifold $M$ with odd complex dimension there may exist an integrable $%
\left( -1,1\right) $-structure (for instance a $\left( 0,1\right) $%
-structure, i.e. a pseudo-K\"{a}hler structure), but there cannot exist $%
\left( -1,0\right) $-structures on $M$ (even non-integrable ones), since
these require $M$ to have even complex dimension. The only other slash
structures that are defective in this sense are $\left( 1,-1\right) $%
-integrable structures on $M$ with odd complex dimension, since $M$ carries
a compatible $\mathbb{C}$-symplectic structure only if its complex dimension
is even.

On the other hand, a stronger possible notion of interpolation on a complex
manifold $M$ could require that, pointwise (or equivalently, at the linear
algebra level on the extended tangent space at a fixed point of $M$), $%
\left( \lambda ,0\right) $- and $\left( 0,\ell \right) $-structures are in
the same orbit of the group $G$ as in Theorem \ref{complexFibre}. The
signature makes this fail for $\left( 1,1\right) $-structures. Indeed, by
that theorem and Proposition \ref{signC}, a pseudo K\"{a}hler structure on $%
M $ is pointwise in the same orbit as a complex product structure (both
compatible with $j$) only if it is split. We have this situation for no
other slash structure;\ in particular, pointwise, hypercomplex and pseudo-K%
\"{a}hler structures on $M$ of any signature (if existing) are in the same $%
G $-orbit.

\subsection{$B$-fields preserving slash structures on $\left( M,j\right) $}

Let $\omega $ be a closed two-form on a $M$ and let $B_{\omega }$ be the
vector bundle isomorphism of $\mathbb{T}M$ defined by $B_{\omega }\left(
u+\sigma \right) =u+\sigma +\omega ^{\flat }\left( u\right) $, which is
called a $B$\emph{-field transformation}. It is well-known that $B_{\omega }$
is an isometry for $b$ and preserves generalized complex and paracomplex
structures (acting by conjugation $S\mapsto B_{\omega }\cdot S=B_{\omega
}\circ S\circ B_{-\omega }$).

\begin{proposition}
Let $\left( M,j\right) $ be a complex manifold and let $\omega $ be a closed
two-form on $M$. If $j$ is symmetric for $\omega $, then $B_{\omega }$
preserves integrable $\left( -1,1\right) $- and $\left( 1,1;n\right) $%
-structures on $M$. Also, if $j$ is skew-symmetric for $\omega $, then $%
B_{\omega }$ preserves integrable $\left( \lambda ,-1\right) $-structures on 
$M$.
\end{proposition}

For instance, a compatible K\"{a}hler form $\omega $ on $\left( M,j\right) $
provides a $B$-field transformation of hypercomplex~/~$\mathbb{C}$%
-symplectic structures on $\left( M,j\right) $, but in general $\omega $
does not need to be nondegenerate..

\medskip

\noindent \textbf{Proof. }Let $\omega $ be as in the statement of the
proposition. To see that $B_{\omega }$ preserves integrable $\left( \lambda
,\ell \right) $-structures on $M$, it suffices to show that $B_{\omega }$
commutes with $J_{\ell }$, or equivalently, that $\omega ^{\flat }j=\ell
j^{*}\omega ^{\flat }$. That is, $j$ is symmetric or skew-symmetric for $%
\omega $, depending on whether $\ell =1$ or $\ell =-1$, which is true by
hypothesis.

Now, let $S$ be an integrable $\left( 1,1;n\right) $-structure on $M$. Since 
$B_{\omega }$ commutes with $J_{+}$ and is an isometry for $b$, one computes 
$\beta _{B_{\omega }\cdot S}=B_{-\omega }^{\ast }\beta _{S}$, and so $\beta
_{B_{\omega }\cdot S}$ and $\beta _{S}$ have the same signature. Thus, $%
B_{\omega }\cdot S$ is an integrable $\left( 1,1;n\right) $-structure.
\hfill $\square $

\subsection{Some examples}

\noindent1) Let $\pi $ be a Poisson structure on a complex manifold $\left(
M,j\right) .$ Then the associated generalized paracomplex structure $S$ on $%
M $ defined by $S\left( u+\sigma \right) =\left( u+\pi ^{\sharp }\sigma
,-\sigma \right) $ (see Example 3 in \cite{Wade}) is not an integrable $%
\left( 1,\ell \right) $-structure for $\ell =\pm 1$, since $S$ does not
anticommute with $J_{\ell }$. (For a bilinear map $\pi :V^{\ast }\times
V^{\ast }\rightarrow \mathbb{R}$, $\pi ^{\sharp }:V^{\ast }\rightarrow V$ is
defined by 
\begin{equation}
\eta \left( \pi ^{\sharp }\left( \xi \right) \right) =\pi \left( \xi ,\eta
\right) \text{,}  \label{numeral}
\end{equation}%
for all $\xi ,\eta \in V^{\ast }$.)

\medskip

\noindent2) Let $M$ be the Lie group $H\times \mathbb{R}$, where $H$ is the
three dimensional Heisenberg group, and let $e$ denote its identity element.
We consider below a left invariant complex structure $j$ on $M$. Not every
left invariant almost symplectic structure on $M$ for which $j$ is
skew-symmetric is integrable, for instance $\theta $ in (\ref{omegae})
below. In particular, given a constant $\left( 1,1;2\right) $-structure on $%
T_{e}M\oplus T_{e}M^{\ast }$, its left invariant extension to $TM\oplus
TM^{\ast }$ is not necessarily integrable.

We exhibit a one-parameter family of integrable $\left( 1,1;2\right) $%
-structures on $M$ such that most of its members are not pure neutral K\"{a}%
hler or complex product structures on $M$, and also they are not obtained
from them via $B$-field transformations.

Let $\mathcal{B}=\left\{ e_{1},e_{2},e_{3},e_{4}\right\} $ be an ordered
basis of Lie~$\left( M\right) $ satisfying $\left[ e_{1},e_{2}\right] =-%
\left[ e_{2},e_{1}\right] =e_{3}$, and the remaining Lie brackets $\left[
e_{i},e_{j}\right] =0$. Let $\mathcal{B}^{\ast }=\left\{
e^{1},e^{2},e^{3},e^{4}\right\} $ be the basis dual to $\mathcal{B}$. It is
easy to see that the left invariant $2$-form $\omega $ on $M$ defined at the
identity by 
\begin{equation}
\theta_{e}=e^{1}\wedge e^{2}+e^{3}\wedge e^{4}  \label{omegae}
\end{equation}%
is not closed. Consider the matrices%
\begin{eqnarray}
J &=&\left( 
\begin{array}{cc}
i & 0 \\ 
0 & i%
\end{array}%
\right) \text{ \ \ \ \ \ and \ \ \ \ \ }R=\left( 
\begin{array}{cc}
r & 0 \\ 
0 & r%
\end{array}%
\right) \text{,}  \label{JotaErre} \\
\text{where\ \ \ \ \ \ \ \ }i &=&\left( 
\begin{array}{cc}
0 & -1 \\ 
1 & 0%
\end{array}%
\right) \text{ \ \ \ \ \ and \ \ \ \ }r=\left( 
\begin{array}{cc}
1 & 0 \\ 
0 & -1%
\end{array}%
\right) \text{.}  \notag
\end{eqnarray}%
Examples 6.4 and 6.5 in \cite{AdrianS} tell us that $J$ and $R$ are the
matrices (with respect to $\mathcal{B}$) of a complex and a paracomplex
structure on $M$, respectively, yielding a complex product structure on $M$
(all of them left invariant). Hence,%
\begin{equation*}
S=\left( 
\begin{array}{cc}
R & 0 \\ 
0 & -R%
\end{array}%
\right)
\end{equation*}%
is the matrix of a left invariant integrable $\left( 1,1;2\right) $%
-structure on $M$, with respect to the oriented basis $\mathcal{C}$ of $%
T_{e}M\oplus T_{e}M^{\ast }$ obtained by juxtaposition of $\mathcal{B}$ with 
$\mathcal{B}^{\ast }$. By Theorem \ref{complexFibre}, a group $G$ isomorphic
to $O\left( 2,\mathbb{C}\right) $ acts by conjugation on the constant $%
\left( 1,1\right) $-structures on $T_{e}M\oplus T_{e}M^{\ast }\cong \mathbb{C%
}^{4}$ (isomorphism determined by $\mathcal{B}$ and $J$). Therefore, if $%
g\in G$, then $gSg^{-1}$ defines a (possibly not integrable) left invariant $%
\left( 1,1;2\right) $-structure on $M$. Let 
\begin{equation*}
D=\left( 
\begin{array}{cc}
0 & d \\ 
d & 0%
\end{array}%
\right) \text{, \ \ \ \ \ \ where \ \ \ \ \ \ }d=\left( 
\begin{array}{cc}
0 & -r \\ 
r & 0%
\end{array}%
\right) \text{.}
\end{equation*}

\begin{proposition}
\label{ex}For any $t\in \mathbb{R}$, $S\left( t\right) =e^{tD}Se^{-tD}$
defines an integrable left invariant $\left( 1,1;2\right) $-structure on $M$%
. If $4t\not\in \mathbb{Z}\pi $, then $S\left( t\right) $ is not trivial as
in the examples in \emph{(\ref{trivialC})} and cannot be obtained from them
via a $B$-field transformation.
\end{proposition}

\noindent \textbf{Proof. }We compute 
\begin{equation*}
e^{tD}=\left( \cos t\right) I_{8}+\left( \sin t\right) D
\end{equation*}%
(we denote by $I_{n}$ the $n\times n$ identity matrix). We then compute%
\begin{equation*}
S\left( t\right) =\left( 
\begin{array}{cc}
\cos \left( 2t\right) R & -\sin \left( 2t\right) T \\ 
\sin \left( 2t\right) T & -\cos \left( 2t\right) R%
\end{array}%
\right) \text{,\ \ \ \ \ \ \ \ \ \ where }\ \ \ \ \ \ \ \ \ \ \ T=\left( 
\begin{array}{cc}
0 & -I_{2} \\ 
I_{2} & 0%
\end{array}%
\right) \text{.}
\end{equation*}%
Then $S\left( t\right) $ is trivial if and only if $4t\in \mathbb{Z}$. In
order to see that (\ref{condition1}) and (\ref{condition2}) are satisfied,
with $\lambda =\ell =1$\ (and also that $S\left( t\right) $ has signature $2$%
), one could check that $D\in $ Lie $\left( G\right) $. We find it simpler
to verify those conditions directly. For this, one uses that the matrices of 
$J_{+}$ and $b$ with respect to the ordered basis $\mathcal{C}$ are%
\begin{equation*}
\left( 
\begin{array}{cc}
J & 0 \\ 
0 & -J%
\end{array}%
\right) \text{\ \ \ \ \ and \ \ \ \ \ \ }\left( 
\begin{array}{cc}
0 & I_{4} \\ 
I_{4} & 0%
\end{array}%
\right) \text{,}
\end{equation*}%
respectively. It would be cumbersome to prove the Courant integrability
condition for $S\left( t\right) $ by definition. Luckily, we can use
Proposition 1.2 in \cite{Izu}: It suffices to check that if

\begin{equation*}
T=\left[ \sigma ^{\flat }\right] _{\mathcal{B},\mathcal{B}^{\ast }}\ \ \ \ \
\ \ \ \ \ \text{and \ \ \ \ \ \ \ \ \ }R=\left[ A\right] _{\mathcal{B},%
\mathcal{B}}\text{,}
\end{equation*}%
then the left invariant extensions of $\sigma $ and $\omega $ to $M$ are
symplectic forms, where $\omega ^{\flat }=\sigma ^{\flat }\circ A$. The
matrix of $\omega $ with respect to $\mathcal{B}$ is $d$ as above. By the
first row of Table 3.3 in \cite{GO}, we have that%
\begin{equation*}
\left( 
\begin{array}{cccc}
0 & c & a & \pm b \\ 
-c & 0 & -b & \pm a \\ 
-a & b & 0 & 0 \\ 
\mp b & \mp a & 0 & 0%
\end{array}%
\right) \text{,}
\end{equation*}%
with $a^{2}+b^{2}\neq 0$, are matrices inducing left invariant symplectic
forms on $M$ (the signs $\pm $ and $\mp $ are allowed, since $\left(
h,s\right) \mapsto \left( h,-s\right) $ is an automorphism of $M$). Now, $T$
and $d$ have this form, hence the left invariant extensions of $\sigma $ and 
$\omega $ to $M$ are symplectic forms. Consequently, $S\left( t\right) $ is
integrable.

Finally, let $a,b$ be $4\times 4$ skew-symmetric matrices, with $\det a\neq
0 $. Let 
\begin{equation*}
Q=\left( 
\begin{array}{cc}
0 & a^{-1} \\ 
a & 0%
\end{array}%
\right) \text{\ \ \ \ \ \ \ and\ \ \ \ \ \ \ \ }B=\exp \left( 
\begin{array}{cc}
0 & 0 \\ 
b & 0%
\end{array}%
\right) \text{,}
\end{equation*}%
and suppose that $4t\not\in \mathbb{Z}\pi $ and%
\begin{equation*}
S\left( t\right) =BQB^{-1}=\left( 
\begin{array}{cc}
-a^{-1}b & a^{-1} \\ 
a-ba^{-1}b & ba^{-1}%
\end{array}%
\right) \text{.}
\end{equation*}%
This implies that $a^{-1}=-a$, $ab=ba$, and so $a=a\left( I_{4}+b^{2}\right) 
$. Hence, $b^{2}=0$ and then $b=0$, since $b$ is skew-symmetric. Similar but
easier computations show that $S\left( t\right) $ can be obtained from a
pure complex product structure via a $B$-field transformation only if $%
S\left( t\right) $ is trivial.

\begin{remark}
We do not know whether there are complex manifolds admitting integrable $%
\left( \lambda ,\ell \right) $-structure but no integrable $\left( \lambda
,0\right) $- or $\left( 0,\ell \right) $-structures.
\end{remark}

\section{Generalized geometric structures on symplectic manifolds\label{S4}}

\subsection{Geometric structures compatible with $\protect\omega $}

Let $\left( M,\omega \right) $ be a symplectic manifold. We consider the
following geometric structures on $M$ compatible with $\omega $.

\bigskip

\noindent \textbf{Integrable (}$1,0$\textbf{)-structure or bi-Lagrangian
foliation of }$\left( M,\omega \right) $ \cite{RB, bil2, bil1}\textbf{.} It
is given by a paracomplex structure $r$ on $M$ which is skew-symmetric for $%
\omega $. Then the leaves of the eigendistributions of $r$ are complementary
Lagrangian submanifolds. This structure is also called \textbf{para-K\"{a}%
hler} \cite{GarciaR, Aleks} or \textbf{K\"{a}hler }$\mathbb{L}$\textbf{%
-manifold }\cite{HarveyL}.

\medskip

\noindent \textbf{Integrable (}$-1,0$\textbf{)-structure or pseudo-K\"{a}%
hler structure on }$\left( M,\omega \right) $. It is given by a complex
structure $j$ on $M$ which is skew-symmetric for $\omega $. If $g$ denotes
the pseudo-Riemannian metric on $M$ given by $g\left( x,y\right) =\omega
\left( jx,y\right) $, then $\left( M,g,j\right) $ is pseudo-K\"{a}hler.%
\textbf{\ }

\medskip

\noindent \textbf{Integrable (}$0,1$\textbf{)-structure or }$\mathbb{L}$%
\textbf{-symplectic structure on }$\left( M,\omega \right) $. It is given by
a symplectic form $\theta $ on $M$ such that the tensor field $A$ given by $%
\theta ^{\flat }=\omega ^{\flat }\circ A$ satisfies $A^{2}=$ id and is
split;\ in particular, $A$ is symmetric for $\omega $. Then $\Omega =\omega
+\varepsilon \theta $ is an $\mathbb{L}$-symplectic structure on $M$ ($TM$
is a vector space over $\mathbb{L}$ via $\left( a+b\varepsilon \right)
u=au+\varepsilon Av$). This structure may be also called a \textbf{%
bi-symplectic foliation on }$\left( M,\omega \right) $. See Proposition \ref%
{EllSympl} below.

\medskip

\noindent \textbf{Integrable (}$0,-1$\textbf{)-structure or }$\mathbb{C}$%
\textbf{-symplectic structure on }$\left( M,\omega \right) $. It is given by
a symplectic form $\theta $ on $M$ such that the tensor field $A$ given by $%
\theta ^{\flat }=\omega ^{\flat }\circ A$ satisfies $A^{2}=-$ id;\ in
particular, $A$ is symmetric for $\omega $. Then $\Omega =\omega -i\theta $
is a $\mathbb{C}$-symplectic structure on $M$.

\medskip

We also have

\medskip

\noindent \textbf{(}$+$\textbf{)-integrable (}$1,0$\textbf{)-structure or
Lagrangian foliation of }$\left( M,\omega \right) $ \textbf{with a
Lagrangian Ehresmann connection.} It is given by a tensor field $r$ of type $%
\left( 1,1\right) $ on $M$ with $r^{2}=\;$id which is skew-symmetric for $%
\omega $, such that the $1$-eigensection $\mathcal{D}_{+}$ of $r$ is an
integrable distribution. Then $M\rightarrow M/_{\mathcal{D}_{+}}$ is a
Lagrangian foliation with $\mathcal{D}_{-}$ (the $\left( -1\right) $%
-eigensection of $r$) a Lagrangian Ehresmann connection.

\bigskip

All these structures compatible with $\omega $ are well-known except
possibly the $\mathbb{L}$-symplectic ones. In the literature we have found
an example in the recent paper \cite{Datta}: If $\sigma _{1}$ and $\sigma
_{2}$ are as in Theorem A in that article, then one can check that $\sigma
_{1}+\varepsilon \sigma _{2}$ is $\mathbb{L}$-symplectic. As it is the case
for $\mathbb{C}$-symplectic structures, if $\left( M,\omega \right) $ admits
an integrable $\mathbb{L}$-symplectic structure, then $\dim M$ is a multiple
of $4$.

\begin{proposition}
\label{EllSympl}Let $\left( M,\omega \right) $ be a symplectic manifold.
Suppose that the closed two-form $\theta $ on $M$ determines an $\mathbb{L}$%
-symplectic structure on $M$. Then, for $\delta =\pm 1$, the $\delta $%
-eigendis\-tri\-butions $D_{\delta }$ of the tensor field $A=(\omega ^{\flat
})^{-1}\circ \theta ^{\flat }$ are integrable and the restriction of $\omega 
$ to the leaves of both foliations is nondegenerate; in particular, the
leaves are symplectic.

Conversely, suppose that $M$ has two complementary foliations $D_{\delta }$ 
\emph{(}$\delta =\pm 1$\emph{) }with equal dimensions and $\left. \omega
\right| _{D_{\delta }\times D_{\delta }}$ is nondegenerate, and define the
tensor field $A$ on $M$ of type $\left( 1,1\right) $ by $\left. A\right|
_{D_{\delta }}=\delta ~$\emph{id}$_{D_{\delta }}$. Then $\theta ^{\flat
}=\omega ^{\flat }\circ A$ determines an $\mathbb{L}$-symplectic structure
on $M$.
\end{proposition}

\noindent \textbf{Proof. }First, we check that $D_{\delta }$ are involutive
for $\delta =\pm 1$. Since $\omega $ is nondegenerate, it suffices to see
that 
\begin{equation}
\omega \left( A\left[ u,v\right] ,z\right) =\delta \omega \left( \left[ u,v%
\right] ,z\right)  \label{omegaTheta}
\end{equation}
for any locally defined vector fields $u,v,z$ on $M$, with $u,v$ local
sections of $D_{\delta }$. Now, the left hand side equals 
\begin{eqnarray*}
\theta \left( \left[ u,v\right] ,z\right) &=&u\theta \left( v,z\right)
-v\theta \left( u,z\right) +z\theta \left( u,v\right) +\theta \left( \left[
u,z\right] ,v\right) -\theta \left( \left[ v,z\right] ,u\right) \\
&=&u\omega \left( Av,z\right) -v\omega \left( Au,z\right) +z\omega \left(
Au,v\right) +\omega \left( \left[ u,z\right] ,Av\right) -\omega \left( \left[
v,z\right] ,Au\right)
\end{eqnarray*}
(we have used that $\theta $ is closed and $A$ is symmetric for $\omega $).
Since $Au=\delta u$, $Av=\delta v$ and $\omega $ is closed, this is the same
as the right hand side of (\ref{omegaTheta}), as desired. Also, one computes
that $\omega \left( D_{+},D_{-}\right) =0$. Hence the form $\omega $
restricted to $D_{\pm }$ is nondegenerate. Similar arguments yield the
converse. \hfill $\square $

\subsection{Slash structures on $\left( M,\protect\omega \right) $}

\begin{definition}
Let $\left( M,\omega \right) $ be a symplectic manifold. For $k=-1,k=1$ let $%
I_{k}$ be the generalized complex, respectively generalized paracomplex,
structure on $M$ given by 
\begin{equation*}
I_{k}=\left( 
\begin{array}{cc}
0 & k\left( \omega ^{\flat }\right) ^{-1} \\ 
\omega ^{\flat } & 0%
\end{array}
\right) \text{.}
\end{equation*}
\end{definition}

\begin{definition}
Let $\left( M,\omega \right) $ be a symplectic manifold. Given $\lambda =\pm
1$ and $\ell =\pm 1$, a generalized complex structure $S$ \emph{(}for $%
\lambda =-1$\emph{)} or a generalized paracomplex structure $S$ \emph{(}for $%
\lambda =1$\emph{)} on $M$ is said to be an \textbf{integrable}\emph{\ }$%
\left( \lambda ,\ell \right) $\textbf{-structure} on $\left( M,\omega
\right) $ if 
\begin{equation}
SI_{\lambda \ell }=I_{\lambda \ell }S\text{ \ \ \ \ \ and \ \ \ \ \ \ }%
I_{\lambda \ell }S\text{ is split.}  \label{conmSympl}
\end{equation}%
The condition of $SI_{\lambda \ell }$ being split is empty if $\ell =-1$,
since $\left( SI_{-\lambda }\right) ^{2}=-$ \emph{id}.

In the same way, a $\left( +\right) $-generalized paracomplex structure $S$
on $M$ is said to be a $\left( +\right) $\textbf{-integrable}\emph{\ }$%
\left( 1,\ell \right) $\textbf{-structure} on $\left( M,\omega \right) $ if $%
SI_{\ell }=I_{\ell }S$ and $SI_{\ell }$ is split.
\end{definition}

\bigskip

We call $\mathcal{S}_{\omega }\left( \lambda ,\ell \right) $ the set of all
integrable $\left( \lambda ,\ell \right) $-structures on $\left( M,\omega
\right) $, and $\mathcal{S}_{\omega }^{+}\left( 1,\ell \right) $ the set of
all $\left( +\right) $-integrable $\left( 1,\ell \right) $-structures.

\begin{example}
If $r$ and $\theta $ are integrable $\left( \lambda ,0\right) $- and $\left(
0,\ell \right) $-structures on $\left( M,\omega \right) $, respectively,
then easy computations show that 
\begin{equation*}
R=\left( 
\begin{array}{cc}
r & 0 \\ 
0 & -r^{\ast }%
\end{array}%
\right) \text{ \ \ \ and\ \ \ \ }Q=\left( 
\begin{array}{cc}
0 & \lambda (\theta ^{\flat })^{-1} \\ 
\theta ^{\flat } & 0%
\end{array}%
\right)
\end{equation*}%
belong to $\mathcal{S}_{\omega }\left( \lambda ,\ell \right) $. We only
comment that $I_{\lambda }Q$ is split since it consists of the blocks $%
\lambda A$ and $\lambda A^{\ast }$, where $A$ is the split tensor field
associated to $\theta $ as in the definition of integrable $\left(
0,1\right) $-structures above. For this, see the end of the proof of Theorem %
\ref{THMsymplectic}.
\end{example}

The following simple theorem justifies the terminology introduced in the
previous subsection and includes the notion of interpolation. See the
comment at the end of the section.

\begin{theorem}
\label{THMsymplectic}Let $\left( M,\omega \right) $ be a symplectic
manifold. For $\lambda =\pm 1,\ell =\pm 1$, integrable $\left( \lambda ,\ell
\right) $-structures on $\left( M,\omega \right) $ interpolate between
integrable $\left( \lambda ,0\right) $- and $\left( 0,\ell \right) $%
-structures on $\left( M,\omega \right) $, that is, if 
\begin{equation*}
R=\left( 
\begin{array}{ll}
r & 0 \\ 
0 & t%
\end{array}
\right) \text{ \ \ \ and\ \ \ \ }Q=\left( 
\begin{array}{cc}
0 & p \\ 
\theta ^{\flat } & 0%
\end{array}
\right)
\end{equation*}
belong to $\mathcal{S}_{\omega }\left( \lambda ,\ell \right) $, then $r$ and 
$\theta $ are integrable $\left( \lambda ,0\right) $- and $\left( 0,\ell
\right) $-structures on $\left( M,\omega \right) $, respectively.

Also, for $\ell =\pm 1$, $\left( +\right) $-integrable $\left( 1,\ell
\right) $-structures interpolate between $\left( +\right) $-integrable $%
\left( 1,0\right) $- and $\left( 0,\ell \right) $-structures on $\left(
M,j\right) $.
\end{theorem}

\noindent \textbf{Proof. }The first paragraph of the proof of Theorem \ref%
{THMcomplex} applies, in particular $t=-r^{\ast }$, $\theta $ is a closed $2$%
-form and $p=\lambda (\theta ^{\flat })^{-1}$, and also $r$ is a complex
structure on $M$ (for $\lambda =-1$) or a tensor field on $M$ of type $%
\left( 1,1\right) $ with $r^{2}=$ id and involutive eigendistributions (for $%
\lambda =1$).

Suppose first that $R$ as above commutes with $I_{\lambda \ell }$. Hence, $%
-r^{\ast }\left( \omega ^{\flat }\right) =\omega ^{\flat }r$, or
equivalently, $\omega \left( u,rv\right) =-\omega \left( ru,v\right) $ for
all vector fields $u,v$. That is, $r$ is skew-symmetric for $\omega $. This
implies that the eigendistributions of $r$ are Lagrangian for $\omega $ and 
so $r$ is split and a paracomplex strucuture.

Now suppose that $QI_{\lambda \ell }=I_{\lambda \ell }Q$. Since $p=\lambda
(\theta ^{\flat })^{-1}$, we have that 
\begin{equation*}
\lambda (\theta ^{\flat })^{-1}\omega ^{\flat }=\lambda \ell (\omega ^{\flat
})^{-1}\theta ^{\flat }\text{.}
\end{equation*}%
Calling $A=(\omega ^{\flat })^{-1}\theta ^{\flat }$, which is a tensor field
of type $\left( 1,1\right) $ on $M$, the expression above yields $%
A^{-1}=\ell A$, or equivalently, $A^{2}=\ell $ id.

Now we verify that $A$ is symmetric for $\omega $, i.e., $\omega \left(
Au,v\right) =\omega \left( u,Av\right) $, or equivalently, $\omega ^{\flat
}\left( Au\right) \left( v\right) =\omega ^{\flat }\left( u\right) \left(
Av\right) $ for all vector fields $u,v$ on $M$. This is the same as $\omega
^{\flat }A=A^{\ast }\omega ^{\flat }$, which is true since $A^{\ast
}=(\theta ^{\flat })^{\ast }(\omega ^{\flat })^{\ast -1}=\left( -1\right)
^{2}\theta ^{\flat }(\omega ^{\flat })^{-1}$ ($\theta $ and $\omega $ are
both skew-symmetric).

It remains only to show that $A$ is split if\ $\ell =1$. By hypothesis, the
matrix $I_{\lambda }Q=\lambda $ diag~$\left( A,A^{\ast }\right) $ is split ($%
A^{-1}=A$). Since the dimensions of the $1$-eigenspaces of $A$ and $A^{\ast
} $ coincide, $A$ must be split.

The last statement is true by the same reasons as in Theorem \ref{THMcomplex}%
. \hfill\ $\square $

\subsection{Slash structures on $\left( M,\protect\omega \right) $ in
classical terms}

\begin{proposition}
\label{fund}An integrable $\left( \lambda ,\ell \right) $-structure $S$ on a
symplectic manifold $\left( M,\omega \right) $ has the form 
\begin{equation}
S=\left( 
\begin{array}{cc}
A & \lambda \ell B\left( \omega ^{\flat }\right) ^{-1} \\ 
\omega ^{\flat }B & -A^{\ast }%
\end{array}%
\right) \text{,}  \label{STensors}
\end{equation}%
where $A$ and $B$ are endomorphisms of $TM$ satisfying 
\begin{equation*}
\lambda A^{2}+\ell B^{2}=\emph{id}\text{, \ \ \ \ \ \ }AB+BA=0\text{,\ \ \ \
\ \ \ \ }\omega ^{\flat }A=-A^{\ast }\omega ^{\flat }
\end{equation*}%
and, for $\ell =1$, that the following matrix \emph{(}which squares to the
identity\emph{)} is split. 
\begin{equation}
\left( 
\begin{array}{cc}
B & A \\ 
\lambda A & B%
\end{array}%
\right) \text{.}  \label{matriz}
\end{equation}
\end{proposition}

\noindent\textbf{Proof. } Since $S$ is a generalized complex (for $\lambda
=-1$) or paracomplex structure (for $\lambda =1$), by \cite{Crainic} (see
also \cite{Izu}) one has 
\begin{equation}
S=\left( 
\begin{array}{cc}
A & \pi ^{\sharp } \\ 
\theta ^{\flat } & -A^{\ast }%
\end{array}%
\right) \text{,}  \label{Sgeneral}
\end{equation}%
where $\theta $ and $\pi $ are skew-symmetric, $\pi ^{\sharp}$ was defined
in (\ref{numeral}), and $A$ satisfies 
\begin{equation}
A^{2}+\pi ^{\sharp }\theta ^{\flat }=\lambda ~\text{id,\ \ \ \ \ \ \ }\theta
^{\flat }A=A^{\ast }\theta ^{\flat }\text{,\ \ \ \ \ and \ \ \ \ }\pi
^{\sharp }A^{\ast }=A\pi ^{\sharp }\text{.}  \label{tres}
\end{equation}%
Now, using the first equation in (\ref{conmSympl}) we have that 
\begin{equation}
\pi ^{\sharp }\omega ^{\flat }=\lambda \ell (\omega ^{\flat })^{-1}\theta
^{\flat }\text{\ \ \ \ \ \ and\ \ \ \ \ \ }\omega ^{\flat }A=-A^{\ast
}\omega ^{\flat }\text{.}  \label{dos}
\end{equation}%
Putting $B=(\omega ^{\flat })^{-1}\theta ^{\flat }$, we have $\pi ^{\sharp
}=\lambda \ell B(\omega ^{\flat })^{-1}$ and so (\ref{STensors}) holds.
Besides, the second equations in (\ref{tres}) and (\ref{dos}) yield $AB+BA=0$%
. Notice that, in particular, $A^{2}+B^{2}=\left( A+B\right) ^{2}$. The last
assertion corresponds to the fact that $SI_{\lambda \ell }$ must be split if 
$\ell =1$, and follows from the fact that an easy computation shows that $%
\lambda \phi ^{-1}SI_{\lambda }\phi $ equals (\ref{matriz}), where $\phi
:TM\oplus TM\rightarrow \mathbb{T}M$ is defined by $\phi \left( u,v\right)
=\left( u,\omega ^{\flat }v\right) $. \hfill\ $\square $

\medskip

M. Crainic obtained in \cite{Crainic} (see also \cite{Izu}) conditions on $%
A,\theta $ and $\pi $ for $S$ as in (\ref{Sgeneral}) to be Courant
integrable. One can deduce conditions on $A$ and $B$ as in (\ref{STensors})
for the integrability of $S$.

\subsection{A signature associated to integrable $\left( -1,1\right) $%
-structures on $\left( M,\protect\omega \right) $}

\begin{proposition}
Let $S$ be an integrable $\left( -1,1\right) $-structure on a symplectic
manifold $\left( M,\omega \right) $ of dimension $2m$. Then the form $\beta
_{S}$ on $\mathbb{T}M$ defined by $\beta _{S}\left( x,y\right) =b\left(
I_{-}Sx,y\right) $ is symmetric and has signature $\left( 4n,\allowbreak
4m-4n\right) $ for some integer $n$ with $0\leq n\leq m$.
\end{proposition}

\noindent \textbf{Proof. }The form $\beta _{S}$ is symmetric since $S$ and $%
I_{-}$ are skew-symmetric for $b$. One has that $\left( I_{-}S\right) ^{2}=$
id. For $\delta =\pm 1$, let $D_{\delta }$ be the $\delta $-eigensection of $%
I_{-}S$. Since $I_{-}S$ is required to be split, $D_{+}$ and $D_{-}$ have
both dimension $2m$.

One computes $b\left( D_{+},D_{-}\right) =0$. For $\delta =\pm 1$ let $%
b^{\delta }:=\left. b\right| _{D_{\delta }\times D_{\delta }}$ and $\beta
^{\delta }:=\left. \beta _{S}\right| _{D_{\delta }\times D_{\delta }}$. By
the orthogonality lemma (2.30 in \cite{Harvey}), $b^{\delta }$ is
nondegenerate. One computes also $b^{\delta }=\delta \beta ^{\delta }$. Now,
since $I_{-}$ is an isometry for $b$ and preserves $D_{\delta }$, $\beta
^{+}=b^{+}$ has signature $\left( 2n,2m-2n\right) $ for some integer $0\leq
n\leq m$. Then $b^{-}$ has signature $\left( 2m-2n,2n\right) $ ($b$ is
split). Therefore, $\beta ^{-}$ has signature $\left( 2n,2m-n\right) $, and
so the signature of $\beta _{S}$ is $\left( 4n,4m-4n\right) $. \hfill \ $%
\square $

\begin{definition}
An integrable $\left( -1,1\right) $-structure $S$ on $\left( M,\omega
\right) $ as above is called an integrable $\left( -1,1;n\right) $%
-structure, and we write \emph{sig}\thinspace $\left( S\right) =n$. If $m=2n$%
, by the next proposition, the $\left( -1,1;n\right) $-structure is called a 
\emph{(}split K\"{a}hler\emph{)}~/~$\mathbb{L}$-symplectic structure on $%
\left( M,\omega \right) $.
\end{definition}

\begin{proposition}
\label{signS}\emph{a) }Let $j$ be an integrable $\left( -1,0\right) $%
-structure on $\left( M,\omega \right) $. Then 
\begin{equation*}
R=\left( 
\begin{array}{cc}
j & 0 \\ 
0 & -j^{*}%
\end{array}
\right)
\end{equation*}
is a $\left( -1,1;n\right) $-structure on $\left( M,\omega \right) $ if and
only if the pseudo-K\"{a}hler metric $g\left( u,v\right) =\omega \left(
ju,v\right) $ on $M$ has signature $\left( 2n,2m-2n\right) $.

\smallskip

\emph{b) }Let $\theta $ be an integrable $\left( 0,1\right) $-structure on $%
\left( M,\omega \right) $. Then 
\begin{equation*}
Q=\left( 
\begin{array}{cc}
0 & -(\theta ^{\flat })^{-1} \\ 
\theta ^{\flat } & 0%
\end{array}%
\right)
\end{equation*}%
is a $\left( -1,1;n\right) $-structure on $\left( M,\omega \right) $ if and
only if $m=2n$.
\end{proposition}

\noindent \textbf{Proof. }a) One computes 
\begin{equation*}
\beta _{R}\left( u+\sigma ,v+\tau \right) =\omega \left( ju,v\right) +\tau
\left( (\omega ^{\flat })^{-1}j^{\ast }\sigma \right) =g\left( u,v\right)
+h\left( \sigma ,\tau \right) \text{,}
\end{equation*}%
where the symmetric form $h$ on $T^{\ast }M$ is defined by the last
equality. Now, 
\begin{equation*}
\left( (\omega ^{\flat })^{\ast }h\right) \left( z,w\right) =h\left( \omega
^{\flat }z,\omega ^{\flat }w\right) =-\omega ^{\flat }\left( w\right) \left(
jz\right) =\omega \left( jz,w\right) =g\left( z,w\right) \text{,}
\end{equation*}%
since for an integrable $\left( -1,0\right) $-structure $j$ on $\left(
M,\omega \right) $, $j$ is skew-symmetric for $\omega $. Therefore, if $\phi
:TM\oplus TM\rightarrow \mathbb{T}M$ is the isomorphism defined at the end
of the proof of Proposition \ref{fund}, then 
\begin{equation*}
\phi ^{\ast }\beta _{R}\left( \left( u,z\right) ,\left( v,w\right) \right)
=g\left( u,v\right) +g\left( z,w\right) \text{. }
\end{equation*}%
This implies the assertion of (a), since $\phi ^{\ast }\beta _{R}$ and $%
\beta _{R}$ have the same signature.

b) As in the definition of integrable $\mathbb{L}$-symplectic structure, we
call $A=(\omega ^{\flat })^{-1}\theta ^{\flat }$. We compute 
\begin{equation*}
\beta _{Q}\left( u+\sigma ,v+\tau \right) =-\tau \left( Au\right) -\sigma
\left( Av\right) \text{.}
\end{equation*}%
We have used that $\theta ^{\flat }(\omega ^{\flat })^{-1}=A^{\ast }$ (since 
$\theta $ and $\omega $ are skew-symmetric) and that $A^{-1}=A$. Since $A$
is split, locally, there exists a basis $\left\{ u_{1},\dots u_{2m}\right\} $
of $TM$ such that $Au_{i}=u_{i}$ for $1\leq i\leq m$ and $Au_{i}=-u_{i}$ for 
$m<i\leq 2m$. Let $\left\{ \alpha _{1},\dots ,\alpha _{2m}\right\} $ be the
dual basis. Analyzing the signs of $\beta _{Q}\left( u_{i}+\alpha
_{i},u_{i}+\alpha _{i}\right) $ and $\beta _{Q}\left( u_{i}-\alpha
_{i},u_{i}-\alpha _{i}\right) $, one concludes that $\beta _{Q}$ is split,
and this yields (b). \hfill\ $\square $

\subsection{The associated homogeneous bundles over $\left( M,\protect\omega %
\right) $}

Now, as we did in the complex case, we work at the algebraic level. We fix $%
p\in M$ and call $\mathbb{E}=\mathbb{T}_{p}M$. By abuse of notation, in the
rest of the section we write $b$ and $I_{k}$ instead of $b_{p}$ and $\left(
I_{k}\right) _{p}$, omitting the subindex $p$.

\begin{theorem}
\label{symplFibre}Let $\left( M,\omega \right) $ be a symplectic manifold of
dimension $2m$. Then, integrable $\left( \lambda ,\ell \right) $- or $\left(
-1,1;n\right) $-structures on $\left( M,\omega \right) $ are smooth sections
of a fiber bundle over $M$ with typical fiber $G/H$, according to the
following table. 
\begin{equation*}
\begin{tabular}{|c|c|c|c|c|}
\hline
$\lambda $ & $\ell $ & \emph{sig} & $G$ & $H$ \\ \hline
$1$ & $1$ & - & $Gl\left( 2m,\mathbb{R}\right) $ & $Gl\left( m,\mathbb{R}%
\right) \times Gl\left( m,\mathbb{R}\right) $ \\ \hline
$1$ & $-1$ & - & $U\left( m,m\right) $ & $Gl\left( m,\mathbb{C}\right) $ \\ 
\hline
$-1$ & $1$ & $n$ & $U\left( m,m\right) $ & $U\left( n,m-n\right) \times
U\left( m-n,n\right) $ \\ \hline
$-1$ & $-1$ & - & $Gl\left( 2m,\mathbb{R}\right) $ & $Gl\left( m,\mathbb{C}%
\right) $ \\ \hline
\end{tabular}%
\end{equation*}
\end{theorem}

Before proving the theorem we introduce some notation and present a
proposition. Let $\sigma \left( \lambda ,\ell \right) $ denote the set of
all $S\in $ End$\,_{\mathbb{R}}\left( \mathbb{E}\right) $ satisfying 
\begin{equation*}
S^{2}=\lambda \,\text{id, }S\text{ is split, skew-symmetric for }b\text{,
and }SI_{\lambda \ell }=I_{\lambda \ell }S\text{ is split}.
\end{equation*}%
Note that $\left( \mathbb{E},I_{k}\right) $ is a vector space over $\mathbb{C%
}$ (respectively, $\mathbb{L}$) for $k=-1$ (respectively, $k=1$). The notion
of $\mathbb{L}$-Hermitian forms \cite{HarveyL} is analogous to the one of $%
\mathbb{C}$-Hermitian forms (see the beginning of Subsection \ref{SubS}).

\begin{proposition}
\label{sesquiSympl}Let $b_{-}:\mathbb{E}\times \mathbb{E}\rightarrow \mathbb{%
C}$ and $b_{+}:\mathbb{E}\times \mathbb{E}\rightarrow \mathbb{L}$ be defined
by 
\begin{equation*}
b_{-}\left( x,y\right) =b\left( x,y\right) -ib\left( x,I_{-}y\right) \text{
\ \ \ \ \ \ and \ \ \ \ \ \ }b_{+}\left( x,y\right) =b\left( x,y\right)
+\varepsilon b\left( x,I_{+}y\right) .
\end{equation*}
Then $b_{-}$ is split $\mathbb{C}$-Hermitian and $b_{+}$ is $\mathbb{L}$%
-Hermitian \emph{(}with respect to $I_{-},I_{+}$, respectively\emph{)}.

\smallskip

Also, if $S\in $ \emph{End}$\,_{\mathbb{R}}\left( \mathbb{E}\right) $
satisfies $S^{2}=\lambda $ \emph{id} and $I_{\lambda \ell }S$ is split, then 
$S\in \sigma \left( \lambda ,\ell \right) $ if and only if 
\begin{equation}
b_{\lambda \ell }\left( Sx,Sy\right) =-\lambda b_{\lambda \ell }\left(
x,y\right)  \label{BLambdaEle}
\end{equation}%
for any $x,y\in \mathbb{E}$.
\end{proposition}

\noindent \textbf{Proof. }We call $\epsilon _{1}=\varepsilon $ and $\epsilon
_{-1}=i$ (in particular, $\epsilon _{k}^{2}=k$). First, for $k=\pm 1$, one
has to show that 
\begin{equation*}
\epsilon _{k}b_{k}\left( x,y\right) =b_{k}\left( x,I_{k}y\right)
=-b_{k}\left( I_{k}x,y\right) \ \ \ \ \ \ \ \ \text{and\ \ \ \ \ \ \ }%
\overline{b_{k}\left( x,y\right) }=b_{k}\left( y,x\right)
\end{equation*}
for all $x,y$. This follows easily from the definitions and the fact that $%
I_{k}$ is skew-symmetric for $b$. Also, $b_{-}$ is split since $b=$ Re~$%
b_{-} $ is split.

Now we prove the second assertion. Suppose first that $S\in \sigma \left(
\lambda ,\ell \right) $. We call $k=\lambda \ell $. Since $S$ commutes with $%
I_{k}$, we compute (using (\ref{symm}) with $T=S$ and $\mu =\lambda $) 
\begin{eqnarray*}
b_{k}\left( Sx,Sy\right) &=&b\left( Sx,Sy\right) +k\epsilon _{k}b\left(
Sx,I_{k}Sy\right) \\
&=&-\lambda b\left( x,y\right) +k\epsilon _{k}b\left( Sx,SI_{k}y\right) \\
&=&-\lambda b\left( x,y\right) -\lambda k\epsilon _{k}b\left( x,I_{k}y\right)
\\
&=&-\lambda b_{k}\left( x,y\right) \text{.}
\end{eqnarray*}%
Conversely, suppose that $S^{2}=\lambda $ id, $SI_{k}$ is split and (\ref%
{BLambdaEle}) holds. By (\ref{symm}) with $T=S$ and $\mu =\lambda $, $S$ is
skew-symmetric for $b=$ Re~$b_{k}$. Now, for $k=\pm 1$, we compute 
\begin{eqnarray*}
b_{k}\left( x,SI_{k}y\right) &=&\lambda b_{k}\left( S^{2}x,SI_{k}y\right)
=\lambda \left( -\lambda \right) b_{k}\left( Sx,I_{k}y\right) =-\epsilon
_{k}b_{k}\left( Sx,y\right) = \\
&=&-\epsilon _{k}\lambda b_{k}\left( Sx,S^{2}y\right) =-\epsilon _{k}\lambda
\left( -\lambda \right) b_{k}\left( x,Sy\right) =b_{k}\left(
x,I_{k}Sy\right) \text{. }
\end{eqnarray*}%
Since $b_{k}$ is nondegenerate, $S$ commutes with $I_{k}$. Therefore, $S\in
\sigma \left( \lambda ,\ell \right) $. $\hfill \square $

\bigskip

\noindent \textbf{Proof of Theorem \ref{symplFibre}. }We follow the same
scheme as in the proof of Theorem \ref{complexFibre}. We suppose first that $%
\lambda \ell =-1$. By the first assertion in Proposition \ref{sesquiSympl},
there exist complex linear coordinates $\left( \phi _{-}\right) ^{-1}=\left(
z,w\right) :\left( \mathbb{E},I_{-}\right) \rightarrow \mathbb{C}^{2m}$ such
that $B_{-}:=\left( \phi _{-}\right) ^{*}b_{-}$ is given by 
\begin{equation*}
B_{-}\left( \left( z,w\right) ,\left( z^{\prime },w^{\prime }\right) \right)
=\overline{z}^{t}w^{\prime }+\overline{w}^{t}z^{\prime }\text{,}
\end{equation*}
which is equivalent to the standard split Hermitian form $\overline{z}%
^{t}z^{\prime }-\overline{w}^{t}w^{\prime }$. Let $\Sigma \left( \lambda
,\ell \right) $ be the subset of End\thinspace $_{\mathbb{C}}\left( \mathbb{C%
}^{2m}\right) $ corresponding to $\sigma \left( \lambda ,\ell \right) $ via
the isomorphism $\phi _{-}$. Clearly $U\left( m,m\right) $ acts on $\Sigma
\left( +,-\right) $ and $\Sigma \left( -,+\right) $ by conjugation.

\medskip

\textbf{Case} $\left( +,-\right) $\textbf{:} Let $S\in $ End$\,_{\mathbb{C}%
}\left( \mathbb{C}^{2m}\right) $ be defined by $S\left( z,w\right) =\left(
z,-w\right) $. Using the second assertion of Proposition \ref{sesquiSympl}
one verifies that $S$ belongs to $\Sigma \left( +,-\right) $ (since $\ell
=-1 $, there is no need to check that $iS$ is split). For $\delta =\pm 1$,
let $V_{\delta }$ be the $\delta $-eingenspace of $S$, that is, 
\begin{equation*}
V_{+}=\left\{ \left( z,0\right) \mid z\in \mathbb{C}^{m}\right\} \text{ \ \
\ \ and \ \ \ \ \ }V_{-}=\left\{ \left( 0,z\right) \mid z\in \mathbb{C}%
^{m}\right\} \text{.}
\end{equation*}%
Given $A\in Gl\left( m,\mathbb{C}\right) $, if $\tilde{A}\left( z,w\right)
=\left( Az,(\overline{A}^{t})^{-1}w\right) $, then $\tilde{A}\in U\left(
m,m\right) $. This provides an inclusion of $Gl\left( m,\mathbb{C}\right) $
into $U\left( m,m\right) $.

Let $H$ be the isotropy subgroup at $S$. For $A\in Gl\left( m,\mathbb{C}%
\right) $, clearly $\tilde{A}$ commutes with $S$ and so $\tilde{A}\in H$.
Besides, if $L\in U\left( m,m\right) $ commutes with $S$, then $L$ preserves 
$V_{+}$ and $V_{-}$. Hence $L\left( z,w\right) =\left( Az,Bw\right) $ for
some $A,B\in Gl\left( m,\mathbb{C}\right) $. Now, $B^{-1}=\overline{A}^{t}$
since $L$ is an isometry for $B_{-}$, and so $L=\tilde{A}$. Therefore $%
H=Gl\left( m,\mathbb{C}\right) $.

The action is transitive: Let $T\in \Sigma \left( +,-\right) $ and for $%
\delta =\pm 1$ let $W_{\delta }$ be the $\delta $-eigenspace of $T$. By (\ref%
{BLambdaEle}), $W_{\delta }$ is isotropic for $B_{-}$. Let $\beta
:W_{+}\rightarrow \left( W_{-}\right) ^{*}$ be given by $\beta \left(
u\right) \left( v\right) =B_{-}\left( \bar{u},v\right) $, which is an
isomorphism of vector spaces over $\mathbb{C}$. Let $u_{1},\dots ,u_{m}$ be
a basis of $W_{+}$ over $\mathbb{C}$ and let $v_{1},\dots ,v_{m}$ be the
basis of $W_{-}$ dual to $\beta \left( \overline{u_{s}}\right) $, $s=1,\dots
,m$. Let $F:\mathbb{C}^{2m}\rightarrow \mathbb{C}^{2m}$ be given by $F\left(
e_{s},0\right) =u_{s}$ and $F\left( 0,e_{s}\right) =v_{s}$. Then $F\in
U\left( m,m\right) $ and $T=FSF^{-1}$. So the action is transitive.

\medskip

\textbf{Case} $\left( -,+;n\right) $\textbf{:} Write $z=\left(
z_{1},z_{2}\right) ,w=\left( w_{1},w_{2}\right) $, with $z_{1},w_{1}\in 
\mathbb{C}^{n},z_{2},w_{2}\in \mathbb{C}^{m-n}$, $0\leq n\leq m$. Let $S\in $
End$\,_{\mathbb{C}}\left( \mathbb{C}^{2m}\right) $ be defined by 
\begin{equation*}
S\left( z_{1},z_{2},w_{1},w_{2}\right) =\left(
-iw_{1},iw_{2},-iz_{1},iz_{2}\right) \text{.}
\end{equation*}
We have that $S^{2}=-$ id and $iS\left( z_{1},z_{2},w_{1},w_{2}\right)
=\left( w_{1},-w_{2},z_{1},-z_{2}\right) $. For $\delta =\pm 1$, the $\delta 
$-eigenspace of $iS$ is 
\begin{equation*}
V_{\delta }=\left\{ \left( z,\delta r\left( z\right) \right) \mid z\in 
\mathbb{C}^{m}\right\} \cong \mathbb{C}^{m}\text{.}
\end{equation*}
where $r\left( z_{1},z_{2}\right) =\left( z_{1},-z_{2}\right) $ for $%
z_{1}\in \mathbb{C}^{n},z_{2}\in \mathbb{C}^{m-n}$. Hence, $iS$ is split.
One computes that $S$ is an isometry for $B_{-}$. Then, the second assertion
of Proposition \ref{sesquiSympl} implies that $S$ belongs to $\Sigma \left(
-,+\right) $. Now, it turns out that 
\begin{equation*}
\text{Re~}B_{-}\left( iS\left( z_{1},z_{2},w_{1},w_{2}\right) ,\left(
z_{1}^{\prime },z_{2}^{\prime },w_{1}^{\prime },w_{2}^{\prime }\right)
\right) =\text{Re~}\left( \overline{w_{1}}^{t}w_{1}^{\prime }-\overline{w_{2}%
}^{t}w_{2}^{\prime }+\overline{z_{1}}^{t}z_{1}^{\prime }-\overline{z_{2}}%
^{t}z_{2}^{\prime }\right) \text{,}
\end{equation*}
which is a real inner product on $\mathbb{C}^{2m}$ of signature $\left(
4n,4m-4n\right) $. Therefore, $S\in \Sigma \left( -,+;n\right) $.

One verifies that $\beta ^{\delta }:=\left. B_{-}\right\vert _{V_{\delta
}\times V_{\delta }}$ is $\mathbb{C}$-Hermitian with Hermitian signature $%
\left( n,m-n\right) $ for $\delta =1$ and $\left( m-n,n\right) $ for $\delta
=-1$. There is an obvious isomorphism $\psi _{\delta }:\mathbb{C}%
^{m}\rightarrow V_{\delta }$, $\psi _{\delta }\left( z\right) =\left(
z,\delta r\left( z\right) \right) $. Given $A\in U\left( n,m-n\right) $ and $%
B\in U\left( m-n,n\right) $, the map $\left( A,B\right) \mapsto \alpha
_{A,B} $ defines an inclusion of $U\left( n,m-n\right) \times U\left(
m-n,n\right) $ into $U\left( m,m\right) $, where $\alpha
_{A_{1},A_{2}}x=\psi _{\delta }A_{\delta }\left( \psi _{\delta }\right)
^{-1}x$ if $x\in V_{\delta }$.

Now suppose that $\alpha $ is in the isotropy subgroup at $S$ of the action
of $U\left( m,m\right) $, or equivalently, that $\alpha $ is in $U\left(
m,m\right) $ and commutes with $S$. Hence, $\alpha $ preserves $V_{\delta }$
for $\delta =\pm 1$. Then, $\alpha $ must have the form $\alpha _{A,B}$ as
above.

It remains to show that the action is transitive. Let $T\in \Sigma \left(
-,+;n\right) $ and for $\delta =\pm 1$ let $W_{\delta }\ $be the $\delta $%
-eigenspace of $iT$ (it is a complex subspace, since it is the $\left(
-\delta i\right) $-eigenspace of $T$). By (\ref{BLambdaEle}), one has that $%
B_{-}\left( W_{+},W_{-}\right) =0$, and so $\gamma ^{\delta }:=\left.
B_{-}\right| _{W_{\delta }\times W_{\delta }}$ is a nondegenerate $\mathbb{C}
$-Hermitian form on $W_{\delta }$. Since $T\in \Sigma \left( -,+;n\right) $, 
$\gamma ^{+}$ and $\gamma ^{-}$ have Hermitian signature $\left(
n,m-n\right) $ and $\left( m-n,n\right) $, respectively. One uses the Basis
Theorem to see that there exists $F\in U\left( m,m\right) $ such that $%
T=FSF^{-1}$. Therefore, $\Sigma \left( -,+;n\right) $ can be identified with 
$U\left( m,m\right) /\left( U\left( n,m-n\right) \times U\left( m-n,n\right)
\right) $, as desired.

\bigskip

Now assume that $\lambda \ell =1.$ By Proposition \ref{sesquiSympl} there
exist Lorentz linear coordinates $\phi _{+}^{-1}:\mathbb{E}\rightarrow 
\mathbb{L}^{2m}$, such that $B_{+}:=\phi _{+}^{\ast }b_{+}$ has the form 
\begin{equation*}
B_{+}\left( Z,Z^{\prime }\right) =\overline{Z}^{t}Z^{\prime }\text{,}
\end{equation*}
where $Z,Z^{\prime }\in \mathbb{L}^{2m}$. Let $\Sigma \left( \lambda ,\ell
\right) $ be the subset of End\thinspace $_{\mathbb{L}}\left( \mathbb{L}%
^{2m}\right) $ corresponding to $\sigma \left( \lambda ,\ell \right) $ via
the isomorphism $\phi _{+}$.

Let $e=\left( 1-\varepsilon \right) /2$, $\overline{e}=\left( 1+\varepsilon
\right) /2$, which are null Lorentz numbers forming a basis of $\mathbb{L}$.
On has $e^{2}=e,\overline{e}^{2}=\overline{e},e\overline{e}=0$ and $%
\varepsilon e=-e,\varepsilon \overline{e}=\overline{e}$.

By \cite{HarveyL} (Section 3), the group $G$ of transformations preserving $%
B_{+}$ (that is, $\mathbb{L}$-unitary transformations) is isomorphic to $%
Gl\left( 2m,\mathbb{R}\right) $; more precisely, any element of $G$ has the
form $\widehat{A}$ for some $A\in Gl\left( 2m,\mathbb{R}\right) $, where 
\begin{equation}
\widehat{A}\left( xe+y\overline{e}\right) =\left( Ax\right) e+\left(
(A^{t})^{-1}y\right) \overline{e}  \label{Atilde}
\end{equation}%
for all $x,y\in \mathbb{R}^{2m}$. Clearly $Gl\left( 2m,\mathbb{R}\right) $
acts by conjugation on $\Sigma \left( +,+\right) $ and $\Sigma \left(
-,-\right) $.

\medskip

\textbf{Case} $\left( +,+\right) $\textbf{:} Let $S\in $ End$\,_{\mathbb{L}%
}\left( \mathbb{L}^{2m}\right) $ be defined by $S\left( xe+y\overline{e}%
\right) =r\left( x\right) e-r\left( y\right) \overline{e}$, where $x,y\in 
\mathbb{R}^{2m}$ and $r\left( x_{1},x_{2}\right) =\left( x_{1},-x_{2}\right) 
$, with $x_{i}\in \mathbb{R}^{m}$ (in particular, $r^{2}=$ id and $r$ is
split). Hence $\varepsilon S\left( xe+y\overline{e}\right) =-r\left(
x\right) e-r\left( y\right) \overline{e}$. Both $S$ and $\varepsilon S$
square to the identity and are split, as required ($I_{+}$ corresponds to
multiplication by $\varepsilon $ in $\mathbb{L}^{2m}$). We compute 
\begin{eqnarray*}
B_{+}\left( S\left( xe+y\overline{e}\right) ,S\left( x^{\prime }e+y^{\prime }%
\overline{e}\right) \right) &=&-\left( r\left( y\right) \right) ^{t}r\left(
x^{\prime }\right) e-\left( r\left( x\right) \right) ^{t}r\left( y^{\prime
}\right) \overline{e} \\
&=&-B_{+}\left( xe+y\overline{e},x^{\prime }e+y^{\prime }\overline{e}\right) 
\text{,}
\end{eqnarray*}%
since $r^{t}r=$ id. Therefore $S\in \Sigma \left( +,+\right) $. The isotropy
subgroup of the action of $Gl\left( 2m,\mathbb{R}\right) $ at $S$ consists
of the maps $\widehat{A}$ as in (\ref{Atilde}), where $A\left(
x_{1},x_{2}\right) =\left( ax_{1},bx_{2}\right) $ for some $a,b\in Gl\left(
m,\mathbb{R}\right) $, hence, it can be identified with $Gl\left( m,\mathbb{R%
}\right) \times Gl\left( m,\mathbb{R}\right) $.

Now, we see that the action is transitive. Let $T\in \Sigma \left(
+,+\right) $. Since $T$ is $\mathbb{L}$-linear, $T\left( xe+y\overline{e}%
\right) =f\left( x\right) e+g\left( y\right) \overline{e}$ for some linear
endomorphisms $f,g$ of $\mathbb{R}^{2m}$. The condition that $T^{2}=$ id
implies that $f^{2}=g^{2}=$ id. Suppose that $f$ and $g$ have signatures $%
\left( k,2m-k\right) $ and $\left( l,2m-l\right) $, respectively. Since both 
$T$ and $\varepsilon T\left( xe+y\overline{e}\right) =-f\left( x\right)
e+g\left( y\right) \overline{e}$ are split by hypothesis, we have that $%
k+l=2m$ and $2m-k+l=2m$. Hence $k=l=m$ and so $f$ and $g$ are split. Then $f$
is conjugate to $r$ in $Gl\left( 2m,\mathbb{R}\right) $, say $f=crc^{-1}$
with $c\in Gl\left( 2m,\mathbb{R}\right) $. Besides, an easy computation
using that $T$ is an anti-isometry for $B_{+}$ yields $g=-\left(
f^{t}\right) ^{-1}$. Therefore $T=CSC^{-1}$ with $C\left( xe+y\overline{e}%
\right) =c\left( x\right) e+\left( c^{t}\right) ^{-1}\left( y\right) 
\overline{e}$, as desired.

\medskip

\textbf{Case} $\left( -,-\right) $\textbf{:} Let $S\in $ End$\,_{\mathbb{L}%
}\left( \mathbb{L}^{2m}\right) $ be defined by $S\left( xe+y\overline{e}%
\right) =j\left( x\right) e+j\left( y\right) \overline{e}$, where $x,y\in 
\mathbb{R}^{2m}$ and $j\left( x_{1},x_{2}\right) =\left( -x_{2},x_{1}\right) 
$, with $x_{i}\in \mathbb{R}^{m}$ (in particular, $j^{2}=-$ id and $j^{t}j=$
id). Computations analogous to those of the case $\left( +,+\right) $ show
that $S\in \Sigma \left( -,-\right) $ and that the isotropy subgroup of the
action of $Gl\left( 2m,\mathbb{R}\right) $ at $S$ consists of the maps $%
\widehat{A}$ as in (\ref{Atilde}), where $A$ commutes with $j$, that is, $%
A\in Gl\left( m,\mathbb{C}\right) $ via the canonical identification of $%
\left( \mathbb{R}^{2m},j\right) $ with $\mathbb{C}^{m}$. Also, transitivity
of the action follows from similar arguments as in the case $\left(
+,+\right) $. \hfill\ $\square $

\bigskip

Finally, we comment on the strength of the notion of interpolation for slash
structures on symplectic manifolds, in analogy with Subsection \ref{weak}
for complex manifolds. Suppose that the dimension of the symplectic manifold 
$M$ is $m=2n$. If $n$ is odd there may exist integrable $\left( \lambda
,\ell \right) $-structures (for instance $\left( \lambda ,0\right) $%
-structures, i.e. pseudo-K\"{a}hler structures or bi-Lagrangian foliations
compatible with $\omega $), but there cannot exist $\left( 0,\ell \right) $%
-structures on $M$ ($\ell =\pm 1$;\ even not integrable ones), since these
require $n$ to be even.

Moreover, by Theorem \ref{symplFibre} and Proposition \ref{signS},
pointwise, a $\left( -1,0\right) $-structure on $M$ (i.e. a pseudo-K\"{a}%
hler structure on $M$ compatible with $\omega $) is in the same $G$-orbit as
a $\left( 0,1\right) $-structure on $M$ ($G$ as in that theorem) only if the
pseudo-K\"{a}hler structure is split. We have this type of shortcoming for
no other slash structure on $\left( M,\omega \right) $;\ in particular,
pointwise, $\mathbb{C}$-symplectic and pseudo-K\"{a}hler structures on $M$
of any signature (if existing) are in the same $G$-orbit.

\bigskip

Most of the structures considered on complex and symplectic manifolds have
been extensively studied. In the bibliography we refer mainly to those which
are less known or have aroused special interest lately.

\bigskip

\noindent \textsc{f}a\textsc{maf}-\textsc{ciem}\newline
\noindent Medina Allende s/n, Ciudad Universitaria, 5000 C\'{o}rdoba,
Argentina\newline
\noindent salvai@famaf.unc.edu.ar\newline

\end{document}